\documentclass[titlepage,draft,12pt]{article} 
\usepackage{amssymb,amsthm} 
\usepackage[a4paper]{geometry}
\usepackage[usenames,dvipsnames]{xcolor}
\usepackage{pstricks,pst-plot} 

\date{}

\newcommand{\ep}{\varepsilon}
\renewcommand{\qed}{{\penalty 10000\mbox{$\quad\Box$}}}
\newcommand{\re}{\mathbb{R}}
\newcommand{\n}{\mathbb{N}}
\newcommand{\cep}{c_{\ep}}
\newcommand{\gep}{g_{\ep}}
\def\dfrac#1#2{\displaystyle{\frac{#1}{#2}}}
\newcommand{\cinfty}{c_{\infty}}
\newcommand{\uinfty}{u_{\infty}}
\newcommand{\lk}{\lambda_{k}}
\newcommand{\wnk}{w_{n,k}}
\newcommand{\Hhat}{\widehat{H}}

\newtheorem{thm}{Theorem}[section]
\newtheorem{thmbibl}{Theorem}

\newtheorem{rmk}[thm]{Remark}
\newtheorem{prop}[thm]{Proposition}
\newtheorem{defn}[thm]{Definition}

\newtheorem{lemma}[thm]{Lemma}

\title{Kirchhoff equations with strong damping}

\author{Marina Ghisi\vspace{1ex}\\ 
{\normalsize Universit\`a degli Studi di Pisa} \\
{\normalsize Dipartimento di Matematica}\\ 
{\normalsize PISA (Italy)}\\
{\normalsize e-mail: \texttt{ghisi@dm.unipi.it}}
\and
Massimo Gobbino\vspace{1ex}\\ 
{\normalsize Universit\`a degli Studi di Pisa} \\
{\normalsize Dipartimento di Matematica}\\ 
{\normalsize PISA (Italy)}\\  
{\normalsize e-mail: \texttt{massimo.gobbino@unipi.it}}
}

\begin{document}
\maketitle
\begin{abstract}
	We consider Kirchhoff equations with strong damping, namely with a
	friction term which depends on a power of the ``elastic''
	operator.  We address local and global existence of solutions in
	two different regimes depending on the exponent in the friction
	term.
	
	When the exponent is greater than 1/2, the dissipation prevails,
	and we obtain global existence in the energy space assuming only
	degenerate hyperbolicity and continuity of the nonlinear term.
	When the exponent is less than 1/2, we assume strict hyperbolicity
	and we consider a phase space depending on the continuity modulus
	of the nonlinear term and on the exponent in the damping.  In this
	phase space we prove local existence, and global existence if
	initial data are small enough.
	
	The regularity we assume both on initial data and on the nonlinear
	term is weaker than in the classical results for Kirchhoff
	equations with standard damping.
	
	Proofs exploit some recent sharp results for the linearized
	equation and suitably defined interpolation spaces.
	
\vspace{6ex}

\noindent{\bf Mathematics Subject Classification 2010 (MSC2010):}
35L70, 35L80, 35L90.


\vspace{6ex}

\end{abstract}

 
\section{Introduction}

Let $H$ be a separable real Hilbert space.  For every $x$ and $y$ in
$H$, let $|x|$ denote the norm of $x$, and let $\langle x,y\rangle$
denote the scalar product of $x$ and $y$.  Let $A$ be a self-adjoint
linear operator on $H$ with dense domain $D(A)$.  We assume that $A$
is nonnegative, namely $\langle Ax,x\rangle\geq 0$ for every $x\in
D(A)$, so that for every $\alpha\geq 0$ the power $A^{\alpha}x$ is
defined provided that $x$ lies in a suitable domain $D(A^{\alpha})$.

We consider the second order evolution equation
\begin{equation}
	u''(t)+2\delta A^{\sigma}u'(t)+
	m\left(|A^{1/2}u(t)|^{2}\right)Au(t)=0
	\label{pbm:k-eqn}
\end{equation}
with initial data
\begin{equation}
	u(0)=u_{0},
	\hspace{3em}
	u'(0)=u_{1},
	\label{pbm:k-data}
\end{equation}
where $\delta>0$ and $\sigma\geq 0$ are real numbers, and
$m:[0,+\infty)\to[0,+\infty)$ is a given function.  Equations of this
type have long been studied in the literature after Kirchhoff's
monograph~\cite{kirchhoff}.  Let us give a brief outline of the
previous results which are closely related to ours.

\paragraph{\textmd{\emph{Non-dissipative Kirchhoff equations}}}

When $\delta=0$, equation (\ref{pbm:k-eqn}) reduces to
\begin{equation}
	u''(t)+m\left(|A^{1/2}u(t)|^{2}\right)Au(t)=0.
	\label{pbm:k-non-diss}
\end{equation}

This is the abstract setting of the hyperbolic partial differential 
equation
$$u_{tt}(x,t)-m\left(\int_{\Omega}|\nabla 
u(x,t)|^{2}\,dx\right)\Delta u(x,t)=0,$$
with suitable boundary conditions in an open set
$\Omega\subseteq\re^{n}$.  When $n=1$ or $n=2$, and $\Omega$ is
bounded, this equation is a possible model for transversal
small-amplitude vibrations of elastic strings or membranes.  The
nonlocal model was derived by Kirchhoff in~\cite{kirchhoff} after some
mathematical simplifications to the full system of (local) equations
of nonlinear elasticity.

From the mathematical point of view, local existence for 
(\ref{pbm:k-non-diss})--(\ref{pbm:k-data}) is known provided that the 
nonlinearity $m(x)$ is Lipschitz continuous and satisfies the 
strict hyperbolicity condition
\begin{equation}
	m(x)\geq\mu_{1}>0
	\quad\quad
	\forall x\geq 0,
	\label{hp:k-sh}
\end{equation}
and initial data $(u_{0},u_{1})$ are in the space $D(A^{3/4})\times
D(A^{1/4})$.  This result was substantially established by Bernstein
in the pioneering paper~\cite{bernstein}, and then refined by many
authors (see~\cite{ap} for a modern version).  Local existence results
extend easily to the so-called mildly degenerate case, namely when the
nonlinear term satisfies just the degenerate hyperbolicity condition
\begin{equation}
	m(x)\geq 0 \quad\quad \forall x\geq 0,
	\label{hp:k-dh}
\end{equation}
but $m(|A^{1/2}u_{0}|^{2})> 0$.  

Global existence has been proved in a multitude of heterogeneous
special situations, such as analytic data
(\cite{bernstein,as,ds:analytic-1,ds:analytic-2}), quasi-analytic data
(\cite{nishihara,gg:qa}), special nonlinearities (\cite{poho-m}),
dispersive equations and small data
(\cite{gh,ds:dispersive,yamazaki:external,matsu-ruzha}), spectral gap
data or spectral gap operators
(\cite{manfrin,hirosawa,gg:global,gg:qa}).

Global existence in Sobolev spaces is no doubt the main open problem
in the theory of Kirchhoff equations.  Local existence in Sobolev
spaces is an open problem as well in at least three situations: when
the nonlinear term is not Lipschitz continuous, when
$m(|A^{1/2}u_{0}|^{2})=0$ (the so-called really degenerate case,
see~\cite{yamazaki:rdg}), and when initial data are just in the energy
space $D(A^{1/2})\times H$ or in any space weaker that
$D(A^{3/4})\times D(A^{1/4})$.  The interested reader is referred to
the survey~\cite{gg:survey} for more details.

From the technical point of view, a key role is played by the
linearization of (\ref{pbm:k-non-diss}), namely equation
\begin{equation}
	u''(t)+c(t)Au(t)=0,
	\label{pbm:lin-non-diss}
\end{equation}
where $c(t)$ is a time-dependent coefficient. It is well-known 
that (\ref{pbm:lin-non-diss}) admits global solutions for all initial 
data in the energy space $D(A^{1/2})\times H$, or more generally in 
$D(A^{\alpha+1/2})\times D(A^{\alpha})$, provided that $c(t)$ is 
Lipschitz continuous and satisfies the strict hyperbolicity 
condition
\begin{equation}
	0<\mu_{1}\leq c(t)\leq\mu_{2}.
	\label{hp:sh}
\end{equation}

In the case of Kirchhoff equations the coefficient 
$c(t)$ is $m(|A^{1/2}u(t)|^{2})$, hence its Lipschitz continuity is 
related to the boundedness of
$$c'(t)=m'\left(|A^{1/2}u(t)|^{2}\right)\cdot
2\langle Au(t),u'(t)\rangle=
2m'\left(|A^{1/2}u(t)|^{2}\right)\cdot\langle 
A^{3/4}u(t),A^{1/4}u'(t)\rangle.$$

This is the point where both the Lipschitz continuity of $m(x)$ and
the choice of the phase space $D(A^{3/4})\times D(A^{1/4})$ come into
play.

When the coefficient $c(t)$ is not Lipschitz continuous, solutions to
(\ref{pbm:lin-non-diss}) can exhibit a severe derivative loss for all
positive times, as shown in the seminal paper~\cite{dgcs}.  In the
sequel we call (DGCS)-phenomenon this instantaneous loss of
regularity.  As a consequence, when $c(t)$ is not Lipschitz continuous
the initial value problem for (\ref{pbm:lin-non-diss}) is well-posed
only for very regular initial data, such as Gevrey or analytic
functions.  We refer to~\cite{dgcs} for the linear theory, and
to~\cite{gg:der-loss} for its application to Kirchhoff equations.

The existence results for linear equations and the (DGCS)-phenomenon
motivate the three main assumptions in the theory of Kirchhoff
equations, namely Lipschitz continuity of the nonlinearity $m(x)$, strict
hyperbolicity, and initial data at least in $D(A^{3/4})\times
D(A^{1/4})$.  These three main assumptions become three
insurmountable barriers when trying to extend the theory.

\paragraph{\textmd{\emph{Kirchhoff equations with standard damping}}}

Equation (\ref{pbm:k-eqn}) with $\delta>0$ and $\sigma=0$ has long 
been studied starting from the 80s, both in the strictly hyperbolic 
case (see~\cite{debrito,yamada}), and in the mildly degenerate case 
(see~\cite{ny} for nonlinear terms of the form $m(x)=x^{\alpha}$ 
and~\cite{gg:diss} for general Lipschitz continuous nonlinear terms). 

As for local solutions, the state of the art is exactly the same as in
the non-dissipative case.  As for global solutions, it was possible
to prove that they exist in Sobolev spaces (such as $D(A)\times
D(A^{1/2})$ and sometimes $D(A^{3/4})\times D(A^{1/4})$) provided that
initial data are small enough in the same space and $m(x)$ is
Lipschitz continuous.  The interested reader is referred to the
survey~\cite{gg:survey-diss} for further details.

In other words, a standard damping seems to be impotent against the
classical three barriers.  As a further evidence, it was recently
shown in~\cite{gg:dgcs-strong} that the linearized equation with
standard dissipation
$$u''(t)+2\delta u'(t)+c(t)Au(t)=0$$
can exhibit the (DGCS)-phenomenon if the coefficient $c(t)$ is not 
Lipschitz continuous.

\paragraph{\textmd{\emph{Equations with strong damping}}}

Mathematical models with strong damping were proposed
in~\cite{CR}, leading to the linear equation
\begin{equation}
	u''(t)+2\delta A^{\sigma}u'(t)+Au(t)=0,
	\label{pbm:lin-constant}
\end{equation}
which afterwards was rigorously analyzed by many authors from
different points of view, and with different choices of $\sigma$.
Here we just quote the first papers~\cite{CT1,CT2,CT3} (see also the
more recent~\cite{HO}) devoted to analyticity properties of the
semigroup associated to (\ref{pbm:lin-constant}), and the
papers~\cite{I2,I3,shibata} by the dispersive school where
(\ref{pbm:lin-constant}) is considered in the concrete case where
$Au=-\Delta u$ in $\re^{n}$ or special classes of unbounded domains.
Finally, equation (\ref{pbm:lin-constant}) was considered
in~\cite{ggh:sd} in full generality, namely for every $\sigma\geq 0$
and every nonnegative self-adjoint operator $A$.

Concerning Kirchhoff equations with strong dissipation, the first
result we are aware of was obtained by Nishihara
in~\cite{nishihara:k-strong}.  He considered the case $\sigma=1$, and
proved global existence for initial data in $D(A)\times D(A^{1/2})$
provided that the nonlinearity is Lipschitz continuous and satisfies
the degenerate hyperbolicity condition (\ref{hp:k-dh}).  Some years
later, Matos and Pereira~\cite{matos-pereira} improved this result by
showing global existence in the energy space $D(A^{1/2})\times H$ with
the same assumptions on $\sigma$ and $m(x)$.  In other words, in the
case $\sigma=1$ Nishihara overcame one of the three classical barriers
(strict hyperbolicity), and then Matos and Pereira overcame one more
barrier (regularity of initial data).  Both results do not extend to
$\sigma<1$, and seem to exploit the Lipschitz continuity of $m(x)$ in
an essential way.  We discuss this technical point at the beginning of
section~\ref{sec:a-priori}.  

As far as we know, references~\cite{nishihara:k-strong}
and~\cite{matos-pereira} represented the state of the art up to now,
and all subsequent papers on the subject were devoted just to
qualitative properties of these solutions, such as decay estimates
(see for example~\cite{nishihara:decay,ono}).

In this paper we consider equation (\ref{pbm:k-eqn}) in the full 
range $\sigma>0$. Two different regimes appear.
\begin{itemize}
	\item When $\sigma>1/2$ (supercritical regime) the dissipation
	prevails in an overwhelming way.  Indeed, in
	Theorem~\ref{thm:sup-K} we prove global existence in the energy
	space assuming only continuity and degenerate hyperbolicity of the
	nonlinear term.

	\item When $0<\sigma\leq 1/2$ (subcritical regime) the dissipation
	competes with the continuity modulus of the nonlinear term.  In
	Theorem~\ref{thm:sub-K} we assume strict hyperbolicity, and we
	prove
	\begin{itemize}
		\item local existence for all initial data in a phase space
		depending on $\sigma$ and on the continuity modulus of $m$.
	
		\item global existence provided that initial data are small in
		the same space.
	\end{itemize}
	
	Just to give an example, if $m(x)$ is $\beta$-H\"{o}lder continuous
	for some $\beta\in(1-2\sigma,1)$, then we obtain local/global
	existence in phase spaces of the form $D(A^{\alpha+1/2})\times
	D(A^{\alpha})$ with 
	\begin{equation}
		4\alpha\beta\geq 1-2\sigma.  
		\label{ex:abs}
	\end{equation}
	
	Note that this condition allows nonlinear terms which are not
	Lipschitz continuous and initial data with $\alpha<1/4$.
\end{itemize} 

\paragraph{\textmd{\emph{Overview of the technique}}}

Theorem~\ref{thm:sup-K} and Theorem~\ref{thm:sub-K} can be proved in
many different ways (here we use a fixed point argument on the
coefficient $c(t)$, but also fixed point arguments on $u(t)$ or
Galerkin approximations could work) after that three main tools have
been developed.

The first tool is linearization.  The results for the non-dissipative
Kirchhoff equation~(\ref{pbm:k-non-diss}) are based on the theory of
the linearized equation (\ref{pbm:lin-non-diss}).  In the same way,
our results for (\ref{pbm:k-eqn}) rely on the theory of the linear
equation
\begin{equation}
	u''(t)+2\delta A^{\sigma}u'(t)+c(t)Au(t)=0,
	\label{pbm:lin-c(t)}
\end{equation}
where $c(t)$ is a time-dependent coefficient.  This theory was
developed in the recent paper~\cite{gg:dgcs-strong}.  The following
two regimes appeared.
\begin{itemize}
	\item  In the supercritical regime $\sigma>1/2$ the dissipation 
	prevails, and (\ref{pbm:lin-c(t)}) generates a semigroup in the 
	energy space $D(A^{1/2})\times H$ provided that $c(t)$ is 
	measurable and satisfies the degenerate hyperbolicity assumption
	\begin{equation}
		0\leq c(t)\leq\mu_{2}.
		\label{hp:dh}
	\end{equation}

	\item In the subcritical regime $0<\sigma\leq 1/2$, if the strict
	hyperbolicity assumption~(\ref{hp:sh}) is satisfied, the
	dissipation competes with the continuity modulus of $c(t)$.  For
	example, if $c(t)$ is $\gamma$-H\"{o}lder continuous for some
	$\gamma\in(0,1)$, then (\ref{pbm:lin-c(t)}) generates a semigroup
	in the energy space if $\gamma>1-2\sigma$, while
	(\ref{pbm:lin-c(t)}) can exhibit the (DGCS)-phenomenon if
	$\gamma<1-2\sigma$. 
\end{itemize}

Let us consider for example the subcritical regime.  In the
non-dissipative case, or in the case with standard dissipation, the
key point was controlling the Lipschitz constant of $c(t)$.  Now the
linear theory tells us that it is enough to control the
$\gamma$-H\"{o}lder constant of $c(t)$ for some $\gamma>1-2\sigma$.
In turn, this constant depends on the continuity moduli of $m(x)$ and
$|A^{1/2}u(t)|^{2}$.

In order to control the latter, we develop our second tool.  We prove
that, when (\ref{pbm:lin-c(t)}) generates a continuous semigroup, the
continuity modulus of $|A^{1/2}u(t)|^{2}$ depends on initial data
only (and not on $c(t)$), according to the following scheme:
\begin{eqnarray*}
	(u_{0},u_{1})\in D(A^{1/2})\times H & \leadsto & 
	|A^{1/2}u(t)|^{2} \mbox{ is continuous},\\
	(u_{0},u_{1})\in D(A^{3/4})\times D(A^{1/4}) & \leadsto & 
	|A^{1/2}u(t)|^{2} \mbox{ is Lipschitz continuous},\\
	(u_{0},u_{1})\in D(A^{\alpha+1/2})\times D(A^{\alpha})
	 & \leadsto & 
	|A^{1/2}u(t)|^{2} \mbox{ is $4\alpha$-H\"{o}lder continuous}.
\end{eqnarray*}

Thus, if the
nonlinear term is $\beta$-H\"{o}lder continuous for some
$\beta\in(0,1)$, and initial data are in $ D(A^{\alpha+1/2})\times
D(A^{\alpha})$ for some $\alpha\in(0,1/4)$, then $c(t)$ is H\"{o}lder
continuous with exponent $4\alpha\beta$, and when
$4\alpha\beta>1-2\sigma$ the linear theory is well-posed.  This is
where condition (\ref{ex:abs}) appears.

The last tool we develop allows to deal with the equality case
$4\alpha\beta=1-2\sigma$.  Roughly speaking, the reason is that an
initial condition which is ``purely'' in $D(A^{\alpha+1/2})\times
D(A^{\alpha})$ does not exist, in the same way as a vector which is
just in $H$, and nothing more, does not exist.  Indeed, any vector in
$H$ lies actually in some better space, and this better space
guarantees a better continuity modulus for $|A^{1/2}u(t)|^{2}$.  We
state this point formally in Proposition~\ref{prop:data-omega}.

\paragraph{\textmd{\emph{Conclusions}}}

Our results are the best possible ones in the
supercritical regime, where all three barriers have been overcome to
the maximum extent.  In the subcritical regime we extended the theory
beyond two of the three barriers, namely to less regular
nonlinearities and less regular initial data.  The width of this
extension depends on $\sigma$ and it is the best possible one
according to the linear theory, because the linearized equation does
admit the (DGCS)-phenomenon beyond the spaces we considered in
this paper.

Our existence results confirm the general paradigm according to which
Kirchhoff equations are well-posed, at least locally in time, whenever
the (DGCS)-phenomenon is excluded.  What happens beyond remains an
open problem, and deep new ideas are likely to be needed
(see~\cite{gg:der-loss} for a partial result).

We also hope that these results could give an indication about the 
regularizing effects one can expect when adding a strong dissipation 
to quasilinear hyperbolic equations.

\paragraph{\textmd{\emph{Structure of the paper}}}

This paper is organized as follows.  In section~\ref{sec:statements}
we introduce the functional setting and we state our main results
concerning local and global solvability for problem
(\ref{pbm:k-eqn})--(\ref{pbm:k-data}).  In section~\ref{sec:linear} we
recall the results we need from the linear theory developed
in~\cite{gg:dgcs-strong}.  In section~\ref{sec:lin-limit} we study how
solutions to the linear problem (\ref{pbm:lin-c(t)}) depend
continuously on the coefficient $c(t)$. In 
section~\ref{sec:interpolation} we introduce our interpolation 
spaces. In section~\ref{sec:a-priori} we prove a priori estimates for 
Kirchhoff equations after discussing what goes wrong in the classical 
ones. Finally, in section~\ref{sec:proofs-main} we prove our main 
results, actually in a stronger form based on our interpolation spaces.

\setcounter{equation}{0}
\section{Notation and statement of main results}\label{sec:statements}

\paragraph{\textmd{\textit{Functional spaces}}}

Let $H$ be a separable Hilbert space.  Let us assume that $H$ admits a
countable complete orthonormal system $\{e_{k}\}_{k\in\n}$ made by
eigenvectors of $A$.  We denote the corresponding eigenvalues by
$\lambda_{k}^{2}$ (with the agreement that $\lk\geq 0$), so that
$Ae_{k}=\lambda_{k}^{2}e_{k}$ for every $k\in\n$.  We do not assume
the monotonicity of the sequence $\{\lk\}$, even if this assumption
would be reasonable in many applications.  Every $u\in H$ can be
written in a unique way in the form $u=\sum_{k=0}^{\infty}u_{k}e_{k}$,
where $u_{k}=\langle u,e_{k}\rangle$ are the Fourier components of
$u$.  In other words, the Hilbert space $H$ can be identified with the
set of sequences $\{u_{k}\}$ of real numbers such that
$\sum_{k=0}^{\infty}u_{k}^{2}<+\infty$.

We stress that this is just a simplifying assumption, with
substantially no loss of generality.  Indeed, according to the
spectral theorem in its general form (see for example Theorem~VIII.4
in~\cite{reed}), one can always identify $H$ with $L^{2}(M,\mu)$ for a
suitable measure space $(M,\mu)$, in such a way that under this
identification the operator $A$ acts as a multiplication operator by
some measurable function $\lambda^{2}(\xi)$.  All definitions and
statements in the sequel can be easily extended to the general setting
just by replacing the sequence $\{\lk\}$ with the function
$\lambda(\xi)$, and the sequence $\{u_{k}\}$ of Fourier components of
$u$ with the element $\widehat{u}(\xi)$ of $L^{2}(M,\mu)$
corresponding to $u$ under the identification of $H$ with
$L^{2}(M,\mu)$.

Powers of the operator $A$ are just defined as 
$A^{\alpha}u:=\sum_{k=0}^{\infty}\lk^{2\alpha}u_{k}e_{k}$, where of 
course
$$u\in D(A^{\alpha})\Longleftrightarrow
\sum_{k=0}^{\infty}(1+\lk^{4\alpha})u_{k}^{2}<+\infty.$$

It is well-known that $D(A^{\alpha})$ is itself a Hilbert space, with 
norm defined by
$$|u|_{D(A^{\alpha})}^{2}:=|u|^{2}+|A^{\alpha}u|^{2}.$$

In concrete applications, $H$ is the space $L^{2}(\Omega)$ for some open set
$\Omega\subseteq\re^{n}$, and $Au=-\Delta u$, with suitable boundary
conditions.  If the boundary of $\Omega$ is regular enough,
$D(A^{\alpha})$ is the set of functions in the Sobolev space
$H^{2\alpha}(\Omega)$ attaining the boundary conditions in a suitable
sense.

\paragraph{\textmd{\textit{Continuity moduli}}}

Throughout this paper we call \emph{continuity modulus} any continuous
function $\omega:[0,+\infty)\to[0,+\infty)$ such that $\omega(0)=0$, 
$\omega(x)>0$ for every $x>0$, and moreover
\begin{equation}
	x\to\omega(x)\mbox{ is a nondecreasing function},
	\label{hp:omega-monot-0}
\end{equation}
\begin{equation}
	x\to\frac{x}{\omega(x)}\mbox{ is a nondecreasing function.}
	\label{hp:omega-monot}
\end{equation}

From (\ref{hp:omega-monot-0}) and (\ref{hp:omega-monot}) one can
easily deduce that every continuity modulus satisfies
\begin{equation}
	\omega(Lx)\leq\max\{1,L\}\omega(x)
	\quad\quad
	\forall x\geq 0,\quad\forall L\geq 0,
	\label{hp:omega-L}
\end{equation}
and that the composition of two continuity moduli is again a 
continuity modulus. We exploit these two properties several times in 
the sequel. 

A function $c:[0,+\infty)\to\re$ is said to be $\omega$-continuous if
\begin{equation}
	|c(a)-c(b)|\leq
	\omega(|a-b|)
	\quad\quad
	\forall a\geq 0,\ \forall b\geq 0.
	\label{hp:ocont}
\end{equation}

More generally, a function $c:X\to\re$ (with $X\subseteq\re$) is said
to be $\omega$-continuous if it satisfies the same inequality for
every $a$ and $b$ in $X$.

\paragraph{\textmd{\textit{Main results}}}

The main results of this paper concern local and global existence for 
problem (\ref{pbm:k-eqn})--(\ref{pbm:k-data}). The first one gives a 
full answer in the case $\sigma>1/2$, which we call supercritical 
dissipation. 

\begin{thm}[Supercritical dissipation]\label{thm:sup-K}
	Let us consider problem (\ref{pbm:k-eqn})--(\ref{pbm:k-data})
	under the following assumptions:
	\begin{itemize}
		\item $A$ is a self-adjoint nonnegative operator on a
		separable Hilbert space $H$,
	
		\item $m:[0,+\infty)\to[0,+\infty)$ is a continuous function,
		
		\item $\sigma>1/2$ and $\delta>0$ are two real numbers,
		
		\item $(u_{0},u_{1})\in D(A^{1/2})\times H$.
	\end{itemize}

	Then the problem admits at least one global solution
	$$u\in C^{0}\left([0,+\infty),D(A^{1/2})\right)\cap
	C^{1}\left([0,+\infty),H\right).$$
\end{thm}

Our second result addresses the case $\sigma\leq 1/2$, which we call 
subcritical dissipation.

\begin{thm}[Subcritical dissipation]\label{thm:sub-K}
	Let us consider problem (\ref{pbm:k-eqn})--(\ref{pbm:k-data})
	under the following assumptions:
	\begin{itemize}
		\item $A$ is a self-adjoint nonnegative operator on a
		separable Hilbert space $H$,
	
		\item $m:[0,+\infty)\to[0,+\infty)$ satisfies the strict 
		hyperbolicity assumption~(\ref{hp:k-sh}) and the 
		$\omega$-continuity assumption (\ref{hp:ocont}) for some 
		continuity modulus $\omega_{m}(x)$,
		
		\item $\sigma\in(0,1/2]$ and $\delta>0$ are two real 
		numbers,
		
		\item there exists $\alpha\in[0,1/4)$ such that
		\begin{equation}
			\limsup_{x\to 0^{+}}
			\frac{[\omega_{m}(x)]^{4\alpha}}{x^{1-2\sigma}}
			<+\infty,
			\label{hp:limsup-alpha}
		\end{equation}
		
		\item $(u_{0},u_{1})\in D(A^{\alpha+1/2})\times 
		D(A^{\alpha})$ for the same $\alpha\in[0,1/4)$ which appears 
		in (\ref{hp:limsup-alpha}).
	\end{itemize}
	
	Then the following conclusions hold true.
	\begin{enumerate}
		\renewcommand{\labelenumi}{(\arabic{enumi})}
		\item  \emph{(Local existence)} There exists $T>0$ such that 
		problem (\ref{pbm:k-eqn})--(\ref{pbm:k-data}) admits at least 
		one local solution
		$$u\in C^{0}\left([0,T],D(A^{\alpha+1/2})\right)\cap
		C^{1}\left([0,T],D(A^{\alpha})\right).$$
	
		\item  \emph{(Global existence)} There exists $\ep_{0}>0$ 
		such that, if initial conditions satisfy
		\begin{equation}
			|u_{1}|^{2}+|A^{1/2}u_{0}|^{2}+
			|A^{\alpha}u_{1}|^{2}+|A^{\alpha+1/2}u_{0}|^{2}
			\leq\ep_{0},
			\label{hp:data-small}
		\end{equation}
		then there exists at least one global solution
		$$u\in C^{0}\left([0,+\infty),D(A^{\alpha+1/2})\right)\cap
		C^{1}\left([0,+\infty),D(A^{\alpha})\right).$$
	\end{enumerate}
\end{thm}

\begin{rmk}
	\begin{em}
		In Theorem~\ref{thm:sub-K}, the $\omega$-continuity assumption
		on $m(x)$ can be easily weakened to local $\omega$-continuity,
		meaning that we can limit ourselves to assume that for every
		$R>0$ the function $m(x)$ has continuity modulus
		$\omega_{R}(x)$ in $[0,R]$, and all these continuity moduli
		satisfy (\ref{hp:limsup-alpha}) with the same exponent $\alpha$
		(of course the limsup may depend on $R$ in an arbitrary way). 
		We spare the reader from this standard generalization.
	\end{em}
\end{rmk}

\begin{rmk}\label{rmk:border-beta}
	\begin{em}
		As a simple example of application of Theorem~\ref{thm:sub-K},
		let us consider the case where $m(x)$ is $\beta$-H\"{o}lder
		continuous, namely $\omega_{m}(x)=Mx^{\beta}$ for some
		constants $M>0$ and $\beta\in(0,1)$.
		
		If $\beta\in(1-2\sigma,1)$, then (\ref{hp:limsup-alpha}) is 
		satisfied if we take $\alpha:=(1-2\sigma)/(4\beta)$. This 
		exponent is less than 1/4, so that we obtain solvability 
		beyond two of the classical three barriers.
		
		If $\beta\in(0,1-2\sigma)$, then there is no $\alpha\leq 1/4$ 
		for which (\ref{hp:limsup-alpha}) is satisfied, so that our 
		theory does not apply to this case.
		
		In the limit case $\beta=1-2\sigma$, one should take
		$\alpha=1/4$ in order to fulfil (\ref{hp:limsup-alpha}).  This
		case in not included in Theorem~\ref{thm:sub-K} above.
		Nevertheless, a careful inspection of the proof reveals that
		in this limit case one can obtain both local and global
		solvability in $D(A^{3/4})\times D(A^{1/4})$ provided that
		initial data are small enough (so that the smallness is 
		required also for local existence). We spare the reader from 
		this subtlety.
	\end{em}
\end{rmk}

\begin{rmk}
	\begin{em}
		The following pictures provide a rough description of the 
		state of the art. In the horizontal axis we represent the 
		regularity of $m(x)$. With some abuse of notation, values 
		$\beta\in(0,1)$ mean that $m(x)$ is $\beta$-H\"{o}lder 
		continuous, $\beta=1$ means that it is Lipschitz continuous, 
		$\beta>1$ means even more regular. The value $\alpha$ in the 
		vertical axis represents the phase space 
		$D(A^{\alpha+1/2})\times D(A^{\alpha})$.

		In the case $\sigma=0$ (not addressed in this paper), both
		local existence and global existence for small data are known
		in the region where $\alpha\geq 1/4$ and $\beta\geq 1$.  
		
		In the case $0<\sigma<1/2$, Theorem~\ref{thm:sub-K} gives the
		same conclusions in a larger region, which is bounded by the
		line $\alpha=(1-2\sigma)/4$, the line $\beta=1-2\sigma$, and
		the curve $\alpha=(1-2\sigma)/(4\beta)$.  The vertical part of
		the boundary is not included, or at least it is just partially
		included as explained at the end of
		Remark~\ref{rmk:border-beta}.
		
		\def\grafico{\psplot{0.315}{0.8}{0.3 x 0.6 exp div}}
		\psset{unit=9ex}
		\noindent
		\pspicture(-0.5,-1)(2.5,3.5)
		\psframe*[linecolor=yellow](0.8,0.6)(2,2)
		\psline[linewidth=0.7\pslinewidth]{->}(-0.3,0)(2.3,0)
		\psline[linewidth=0.7\pslinewidth]{->}(0,-0.3)(0,2.3)
		\psline[linewidth=0.7\pslinewidth,linestyle=dashed](0.8,0)(0.8,2)
		\psline[linewidth=0.7\pslinewidth,linestyle=dashed](0,0.6)(2,0.6)
		\psline[linewidth=1.5\pslinewidth](0.8,2)(0.8,0.6)(2,0.6)
		\uput[-90](2.1,0){$\beta$}
		\uput[180](0,2.1){$\alpha$}
		\uput[-90](0.8,0){$1$}
		\uput[180](0,0.6){$\frac{1}{4}$}
		\rput(1,-0.6){$\sigma=0$}
		\psframe*[linecolor=yellow](-0.3,3)(-0.05,3.25)
		\rput[Bl](0.1,3){Strict hyperbolicity, local existence for 
		all data, global 
		existence for small data}
		\psframe*[linecolor=green](-0.3,2.6)(-0.05,2.85)
		\rput[Bl](0.1,2.6){Degenerate hyperbolicity, global existence for all data}
		\endpspicture
		\hfill
		
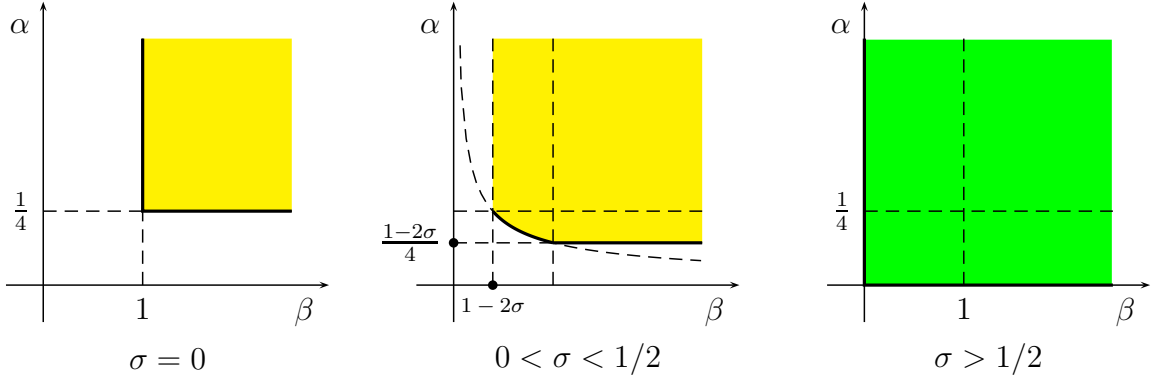
\pspicture(-0.5,-1)(2.5,3.5)
		\pscustom[fillstyle=solid,fillcolor=yellow,linestyle=none]{
		\psline(0.315,2)(0.315,0.6)
		\grafico(2,0.6)(2,2)
		\psline(0.8,0.343)(2,0.343)(2,2)}
		{\psset{linewidth=1.5\pslinewidth}\grafico\psline(0.8,0.343)(2,0.343)}%
		\psclip{\psframe[linestyle=none](0,0)(2,2)}
		\psplot[linewidth=0.7\pslinewidth,linestyle=dashed]{0.02}{2}{0.3 x 0.6 exp div}
		\endpsclip
		\psline[linewidth=0.7\pslinewidth]{->}(-0.3,0)(2.3,0)
		\psline[linewidth=0.7\pslinewidth]{->}(0,-0.3)(0,2.3)
		\psline[linewidth=0.7\pslinewidth,linestyle=dashed](0,0.343)(2,0.343)
		\psline[linewidth=0.7\pslinewidth,linestyle=dashed](0.315,0)(0.315,2)
		\psline[linewidth=0.7\pslinewidth,linestyle=dashed](0.8,0)(0.8,2)
		\psline[linewidth=0.7\pslinewidth,linestyle=dashed](0,0.6)(2,0.6)
		\psdots(0.315,0)(0,0.343)
		\uput[-90](2.1,0){$\beta$}
		\uput[180](0,2.1){$\alpha$}
		\uput[-90](0.315,0){{\scriptsize $1-2\sigma$}}
		\uput[180](0,0.343){$\frac{1-2\sigma}{4}$}
		\rput(1,-0.6){$0<\sigma<1/2$}
		\endpspicture
		\hfill
		\pspicture(-0.5,-1)(2.5,3.5)
		\psframe[fillstyle=solid,fillcolor=green,linestyle=none](0,0)(2,2)
		\psline[linewidth=0.7\pslinewidth]{->}(-0.3,0)(2.3,0)
		\psline[linewidth=0.7\pslinewidth]{->}(0,-0.3)(0,2.3)
		\psline[linewidth=0.7\pslinewidth,linestyle=dashed](0.8,0)(0.8,2)
		\psline[linewidth=0.7\pslinewidth,linestyle=dashed](0,0.6)(2,0.6)
		\psline[linewidth=1.5\pslinewidth](2,0)(0,0)(0,2)
		\uput[-90](2.1,0){$\beta$}
		\uput[180](0,2.1){$\alpha$}
		\uput[-90](0.8,0){$1$}
		\uput[180](0,0.6){$\frac{1}{4}$}
		\rput(1,-0.6){$\sigma>1/2$}
		\endpspicture
	\end{em}
\end{rmk}

In the case $\sigma>1/2$, Theorem~\ref{thm:sup-K} guarantees
global existence in the phase space, independently of the size
of initial data, even in the degenerate hyperbolic case.

\begin{rmk}
	\begin{em}
		Let us discuss in detail the borderline case $\sigma=1/2$.
		First of all, it fits in the framework of
		Theorem~\ref{thm:sub-K}, where for $\sigma=1/2$ assumption
		(\ref{hp:limsup-alpha}) holds true for free, even if
		$\alpha=0$.  Thus we obtain both local existence and global
		existence for small data in the energy space $D(A^{1/2})\times
		H$, provided that $m(x)$ is continuous and satisfies the
		strict hyperbolicity condition (\ref{hp:k-sh}).
		
		Concerning the degenerate hyperbolic case, the limit exponent
		$\sigma=1/2$ is not included in Theorem~\ref{thm:sup-K} as
		stated above.  Nevertheless, a careful inspection of the proof
		reveals that some conclusions hold true also when
		$\sigma=1/2$.  Indeed, we can prove that a local solution
		exists for all initial data $(u_{0},u_{1})$ in the energy
		space such that $4\delta^{2}>m(|A^{1/2}u_{0}|^{2})$, and this
		solution survives as long as
		$4\delta^{2}>m(|A^{1/2}u(t)|^{2})$.  This could be always the
		case for suitable choices of $m(x)$ and $(u_{0},u_{1})$, for
		example in the trivial case where $m(x)<4\delta^{2}$ for every
		$x\geq 0$.  We spare the reader from this further subtlety.
		
		This discussion shows that in the critical case $\sigma=1/2$
		we have coexistence of two statements, one dealing with
		strictly hyperbolic nonlinearities and one dealing with
		degenerate hyperbolic nonlinearities.  This coexistence
		reflects the analogous coexistence in the linear case
		(see~\cite[Remark~3.7]{gg:dgcs-strong}) and suggests the
		existence of a more general unifying statement, which could
		probably deserve further investigation.
	\end{em}
\end{rmk}

\begin{rmk}
	\begin{em}
		Solutions to (\ref{pbm:k-eqn}) are also solutions to
		(\ref{pbm:lin-c(t)}), hence they inherit all properties of
		solutions to linear equations.  In particular, multiple
		choices of the phase space are possible when $\sigma>1/2$, as
		observed in~\cite{ggh:sd,gg:dgcs-strong}.  Therefore, under
		the assumptions of Theorem~\ref{thm:sup-K}, the problem is
		well-posed not only in the energy space or more generally in
		$D(A^{\alpha+1/2})\times D(A^{\alpha})$, but also in all phase
		spaces of the form $D(A^{\alpha})\times D(A^{\beta})$ with
		$\alpha\geq 1/2$, $\beta\geq 0$, and
		$1-\sigma\leq\alpha-\beta\leq\sigma$.
	\end{em}
\end{rmk}

\setcounter{equation}{0}
\section{The linearized equation -- Previous results}\label{sec:linear}

In this section we collect all the results we need from the linear
theory developed in~\cite{gg:dgcs-strong}.  To this end, we introduce
some further notation.  Given any $\nu\geq 0$, we write $H$ as an
orthogonal direct sum
\begin{equation}
	H:=H_{\nu,-}\oplus H_{\nu,+},
	\label{defn:H+-}
\end{equation}
where $H_{\nu,-}$ is the closure of the subspace generated by all
eigenvectors of $A$ relative to eigenvalues $\lk<\nu$, and $H_{\nu,+}$
is the closure of the subspace generated by all eigenvectors of $A$
relative to eigenvalues $\lk\geq\nu$.  For every vector $u\in H$, we
write $u_{\nu,-}$ and $u_{\nu,+}$ to denote its components with
respect to the decomposition (\ref{defn:H+-}).  We point out that
$H_{\nu,-}$ and $H_{\nu,+}$ are $A$-invariant subspaces of $H$, and
that $A$ is a bounded operator when restricted to $H_{\nu,-}$, and a
coercive operator when restricted to $H_{\nu,+}$ if $\nu>0$.

Here we state the results in local form, namely with coefficients 
$c(t)$ defined in some time-interval $[0,T]$. They follow immediately 
from the corresponding results proved in~\cite{gg:dgcs-strong} for 
coefficients defined in $[0,+\infty)$ (it is enough to extend the 
coefficient by setting $c(t)=c(T)$ for every $t\geq T$). 

\subsection{Existence and estimates in the energy space}

As always, we distinguish the supercritical and the subcritical case. 
In both cases, we limit ourselves to initial data in the energy space.

\begin{thmbibl}[Supercritical dissipation,
	see~{\cite[Theorem~3.1]{gg:dgcs-strong}}]\label{thmbibl:sup} 
	Let $T>0$, and let us consider problem
	(\ref{pbm:lin-c(t)})--(\ref{pbm:k-data}) under the following
	assumptions:
	\begin{itemize}
		\item $A$ is a self-adjoint nonnegative operator on a
		separable Hilbert space $H$,
	
		\item the coefficient $c:[0,T]\to\re$ is measurable and
		satisfies the degenerate hyperbolicity assumption
		(\ref{hp:dh}),
	
		\item $\sigma$ and $\delta$ are two positive real numbers such
		that either $\sigma>1/2$, or $\sigma=1/2$ and
		$4\delta^{2}\geq\mu_{2}$,
		
		\item $(u_{0},u_{1})\in D(A^{1/2})\times H$.
	\end{itemize}

	Then the problem has a unique solution 
	\begin{equation}
		u\in C^{0}\left([0,T],D(A^{1/2})\right) \cap
		C^{1}\left([0,T],H\right).
		\label{th:u-reg}
	\end{equation}
	
	Moreover, there exists a constant $K_{1}(\delta,\mu_{2})$,
	depending only on $\delta$ and $\mu_{2}$ (and in particular
	independent of $T$), with the following property.  For every
	$\nu\geq 1$ such that
	\begin{equation}
		4\delta^{2}\nu^{4\sigma-2}\geq\mu_{2},
		\label{hp:sup-nu}
	\end{equation}
	it turns out that
	\begin{equation}
		|u'(t)|^{2}+|A^{1/2}u(t)|^{2}\leq K_{1}(\delta,\mu_{2})
		e^{\nu(1+\mu_{2})t}
		\left(|u_{1}|^{2}+ |A^{1/2}u_{0}|^{2}\right)
		\quad\quad
		\forall t\in [0,T],
		\label{th:sup-sumup}
	\end{equation}
	and more precisely for every $t\in [0,T]$ it turns out that
	\begin{equation}
		|u_{\nu,-}'(t)|^{2}+|A^{1/2}u_{\nu,-}(t)|^{2}\leq 
		e^{\nu(1+\mu_{2})t}\left(|u_{1,\nu,-}|^{2}+
		|A^{1/2}u_{0,\nu,-}|^{2}\right),
		\label{th:sup-u-}
	\end{equation}
	\begin{equation}
		|u_{\nu,+}'(t)|^{2}+|A^{1/2}u_{\nu,+}(t)|^{2}\leq 
		K_{1}(\delta,\mu_{2})\left(|u_{1,\nu,+}|^{2}+
		|A^{1/2}u_{0,\nu,+}|^{2}\right).
		\label{th:sup-u+}
	\end{equation}
\end{thmbibl}

\begin{thmbibl}[Subcritical dissipation, 
	see~{\cite[Theorem~3.2]{gg:dgcs-strong}}]\label{thmbibl:sub}
	Let $T>0$, and let us consider problem
	(\ref{pbm:lin-c(t)})--(\ref{pbm:k-data}) under the following
	assumptions:
	\begin{itemize}
		\item $A$ is a self-adjoint nonnegative operator on a
		separable Hilbert space $H$,
	
		\item  the coefficient $c:[0,T]\to\re$ satisfies the 
		strict hyperbolicity assumption (\ref{hp:sh}) and the 
		$\omega$-continuity assumption (\ref{hp:ocont}) for some 
		continuity modulus $\omega(x)$,
	
		\item $\sigma\in[0,1/2]$ and $\delta>0$ are two real 
		numbers such that
		\begin{equation}
			4\delta^{2}\mu_{1}>\Lambda_{\infty}^{2}+2\delta\Lambda_{\infty},
			\label{hp:sub}
		\end{equation}
		where we set
		\begin{equation}
			\Lambda_{\infty}:=
			\limsup_{\ep\to 0^{+}}\frac{\omega(\ep)}{\ep^{1-2\sigma}},
			\label{defn:Linfty}
		\end{equation}
	
		\item $(u_{0},u_{1})\in D(A^{1/2})\times H$.

	\end{itemize}
	
	Then the problem has a unique solution with the regularity stated 
	in (\ref{th:u-reg}).
			
	Moreover, there exists a constant $K_{2}(\delta,\mu_{1},\mu_{2})$,
	depending only on $\delta$, $\mu_{1}$ and $\mu_{2}$ (and in
	particular independent of $T$), with the following property.  For
	every $\nu\geq 1$ such that
	\begin{equation}
		4\delta^{2}\mu_{1}\geq
		\left[\lambda^{1-2\sigma}\omega\left(\frac{1}{\lambda}\right)\right]^{2}
		+2\delta
		\left[\lambda^{1-2\sigma}\omega\left(\frac{1}{\lambda}\right)\right]
		\quad\quad
		\forall\lambda\geq\nu,
		\label{hp:sub-lambda}
	\end{equation}
	it turns out that
	$$|u'(t)|^{2}+|A^{1/2}u(t)|^{2}\leq K_{2}(\delta,\mu_{1},\mu_{2})
	e^{\nu(1+\mu_{2})t}
	\left(|u_{1}|^{2}+ |A^{1/2}u_{0}|^{2}\right)
	\quad\quad
	\forall t\in [0,T],$$
	and more precisely for every $t\in [0,T]$ it turns out that
	\begin{equation}
		|u_{\nu,-}'(t)|^{2}+|A^{1/2}u_{\nu,-}(t)|^{2}\leq 
		e^{\nu(1+\mu_{2})t}\left(|u_{1,\nu,-}|^{2}+
		|A^{1/2}u_{0,\nu,-}|^{2}\right),
		\label{th:sub-u-}
	\end{equation}
	\begin{equation}
		|u_{\nu,+}'(t)|^{2}+|A^{1/2}u_{\nu,+}(t)|^{2}\leq 
		K_{2}(\delta,\mu_{1},\mu_{2})\left(|u_{1,\nu,+}|^{2}+
		|A^{1/2}u_{0,\nu,+}|^{2}\right).
		\label{th:sub-u+}
	\end{equation}
\end{thmbibl}

The key tool in the proof of Theorem~\ref{thmbibl:sup} and
Theorem~\ref{thmbibl:sub} are some estimates on the Fourier components
of the solution.  We recall these estimates because we need them in
the sequel.  Let $\{u_{k}(t)\}$ denote the components of $u(t)$ with
respect to the usual orthonormal system $\{e_{k}\}$, and let
$\{u_{0k}\}$ and $\{u_{1k}\}$ denote the corresponding components of
initial conditions.  It turns out that $u_{k}(t)$ is a solution to
$$u_{k}''(t)+2\delta\lk^{2\sigma}u_{k}'(t)+
\lk^{2}c(t)u_{k}(t)=0,$$
with initial data
$$u_{k}(0)=u_{0k},
\hspace{3em}
u_{k}'(0)=u_{1k}.$$

Estimates on low-frequency components are quite general 
(see~\cite[Remark~3.3]{gg:dgcs-strong}). Let us assume only that 
$c(t)$ is a measurable function satisfying the degenerate 
hyperbolicity assumption (\ref{hp:dh}) in $[0,T]$, and that 
$\delta\geq 0$, $\sigma\geq 0$, $\nu\geq 0$ and $\lambda_{k}\leq\nu$. 
Then it turns out that
\begin{equation}
	|u_{k}'(t)|^{2}+\lk^{2}|u_{k}(t)|^{2}\leq
	e^{\nu(1+\mu_{2})t}
	\left(|u_{1k}|^{2}+\lk^{2}|u_{0k}|^{2}\right)
	\quad\quad
	\forall t\in[0,T].
	\label{est:uk-}
\end{equation}

Summing over all indices $k$ with $\lambda_{k}<\nu$ we obtain both 
(\ref{th:sup-u-}) and (\ref{th:sub-u-}).

Estimates on high-frequency components are more delicate, and
different in the two cases.  In the supercritical case, let us assume
that $c(t)$, $\delta$ and $\sigma$ are as in
Theorem~\ref{thmbibl:sup}, that $\nu\geq 1$ satisfies
(\ref{hp:sup-nu}), and that $\lambda_{k}\geq\nu$. Then it turns out 
that
\begin{equation}
	|u_{k}'(t)|^{2}+\lk^{2}|u_{k}(t)|^{2}\leq
	K_{1}(\delta,\mu_{2})
	\left(|u_{1k}|^{2}+\lk^{2}|u_{0k}|^{2}\right)
	\quad\quad
	\forall t\in[0,T],
	\label{est:uk+sup}
\end{equation}
where the constant $K_{1}(\delta,\mu_{2})$ depends only on $\delta$ and $\mu_{2}$.
This was established in~\cite[Lemma~5.1]{gg:dgcs-strong}.  Since the
constant is independent of $k$, summing over all indices $k$ with
$\lambda_{k}\geq\nu$ we obtain (\ref{th:sup-u+}).  Finally,
(\ref{th:sup-sumup}) follows from (\ref{th:sup-u-}) and
(\ref{th:sup-u+}) (in this point we need that
$K_{1}(\delta,\mu_{2})\geq 1$, which in turn is a consequence of
(\ref{th:sup-u+}) with $t=0$).

In the subcritical case, let us assume
that $c(t)$, $\delta$ and $\sigma$ are as in
Theorem~\ref{thmbibl:sub}, that $\nu\geq 1$ satisfies
(\ref{hp:sub-lambda}), and that $\lambda_{k}\geq\nu$. Then it turns out 
that
\begin{equation}
	|u_{k}'(t)|^{2}+\lk^{2}|u_{k}(t)|^{2}\leq
	K_{2}(\delta,\mu_{1},\mu_{2})
	\left(|u_{1k}|^{2}+\lk^{2}|u_{0k}|^{2}\right)
	\quad\quad
	\forall t\in[0,T],
	\label{est:uk+sub}
\end{equation}
where the constant $K_{2}(\delta,\mu_{1},\mu_{2})$ depends only on
$\delta$, $\mu_{1}$ and $\mu_{2}$.  This was established
in~\cite[Lemma~5.2]{gg:dgcs-strong}.  Since the constant is
independent of $k$, summing over all indices $k$ with
$\lambda_{k}\geq\nu$ we obtain (\ref{th:sub-u+}).

\subsection{The regularizing effect}\label{sec:regularizing}

The strong dissipation has a regularizing effect in the range
$\sigma\in(0,1)$.  This effect was studied in~\cite{gg:dgcs-strong} in
terms of Gevrey spaces.  Here we just state what we need in the
sequel, without relying on the theory of abstract Gevrey spaces.

Let us start with the supercritical case.
\begin{thmbibl}[Supercritical dissipation -- Regularizing 
	effect]\label{thmbibl:sup-gevrey}
	Let $u(t)$ be a solution to problem 
	(\ref{pbm:lin-c(t)})--(\ref{pbm:k-data}) under the same 
	assumptions of Theorem~\ref{thmbibl:sup}. Let us assume in 
	addition that either $\sigma\in(1/2,1)$, or $\sigma=1/2$ and 
	$4\delta^{2}>\mu_{2}$. 
	
	Let us consider the function
	$$C(t):=\int_{0}^{t}c(s)\,ds
	\quad\quad
	\forall t\in[0,T],$$
	and let us distinguish three cases.
	\begin{enumerate}
		\renewcommand{\labelenumi}{(\arabic{enumi})}
		\item  If $C(t)>0$ for every $t\in(0,T]$, then
		\begin{equation}
			u\in C^{1}\left((0,T],D(A^{\alpha})\right)
			\quad\quad
			\forall\alpha\geq 0.
			\label{th:sup-gevrey-1}
		\end{equation}
	
		\item  If $C(t)=0$ for every $t\in[0,T]$, then
		\begin{equation}
			u'\in C^{\infty}\left((0,T],D(A^{\alpha})\right)
			\quad\quad
			\forall\alpha\geq 0.
			\label{th:sup-gevrey-2}
		\end{equation}
	
		\item If there exists $S_{0}\in(0,T)$ such that $C(t)=0$ for 
		every $t\in[0,S_{0}]$, and $C(t)>0$ for 
		every $t\in(S_{0},T]$, then 
		$$u\in C^{1}\left((S_{0},T],D(A^{\alpha})\right)
			\quad\quad
			\forall\alpha\geq 0,$$
		and
		$$u'\in C^{0}\left((0,S_{0})\cup(S_{0},T],
			D(A^{\alpha})\right)
			\quad\quad
			\forall\alpha\geq 0.$$
	\end{enumerate}
	
\end{thmbibl}

In the first case the regularity (\ref{th:sup-gevrey-1}) follows
from~\cite[Theorem~3.8]{gg:dgcs-strong}.

In the second case it turns out that also $c(t)=0$ for every 
$t\in[0,T]$, so that $u'(t)$ is the solution to the
parabolic problem 
$$v'(t)+2\delta A^{\sigma}v(t)=0, 
\quad\quad
v(0)=u_{1},$$
and therefore the regularity (\ref{th:sup-gevrey-2}) follows from the
regularizing effect of parabolic problems.

In the third case the regularity in $(S_{0},T]$ follows
from~\cite[Theorem~3.8]{gg:dgcs-strong} as in the first case, while
the regularity of $u'(t)$ in $(0,S_{0})$ follows from the parabolic
problem as in the second case (indeed $c(t)=0$ for every
$t\in[0,S_{0}]$).

The last result concerns the regularizing effect in the subcritical
case, and it is an immediate consequence
of~\cite[Theorem~3.9]{gg:dgcs-strong}.

\begin{thmbibl}[Subcritical dissipation -- Regularizing 
	effect]\label{thmbibl:sub-gevrey}
	Let $u(t)$ be a solution to problem 
	(\ref{pbm:lin-c(t)})--(\ref{pbm:k-data}) under the same 
	assumptions of Theorem~\ref{thmbibl:sub}. Let us assume in 
	addition that $\sigma\in(0,1/2]$, so that the case $\sigma=0$ is 
	excluded.
	
	Then (\ref{th:sup-gevrey-1}) holds true.
\end{thmbibl}

\setcounter{equation}{0}
\section{Passing to the limit in linear problems}\label{sec:lin-limit}

In this section we consider a sequence of linear problems
\begin{equation}
	u_{n}''(t)+2\delta A^{\sigma}u_{n}'(t)+c_{n}(t)Au_{n}(t)=0,
	\label{pbm:eqn-cn}
\end{equation}
with fixed initial data
\begin{equation}
	u_{n}(0)=u_{0},
	\hspace{3em}
	u_{n}'(0)=u_{1}.
	\label{pbm:data-cn}
\end{equation}

We assume that the sequence of coefficients
$c_{n}:[0,T]\to[0,+\infty)$ pointwise converges to a limit coefficient
$\cinfty:[0,T]\to[0,+\infty)$, namely
\begin{equation}
	\lim_{n\to +\infty}c_{n}(t)=\cinfty(t)
	\quad\quad
	\forall t\in[0,T].
	\label{hp:cn-cinfty}
\end{equation}

We investigate the convergence of solutions to 
(\ref{pbm:eqn-cn})--(\ref{pbm:data-cn}) to the solution of the limit 
problem
\begin{equation}
	u_{\infty}''(t)+2\delta A^{\sigma}u_{\infty}'(t)+\cinfty(t)Au_{\infty}(t)=0,
	\label{pbm:eqn-cinfty}
\end{equation}
\begin{equation}
	u_{\infty}(0)=u_{0},
	\hspace{3em}
	u_{\infty}'(0)=u_{1}.
	\label{pbm:data-cinfty}
\end{equation}

Our goal is proving uniform convergence of solutions in the energy 
space, namely
\begin{equation}
	\lim_{n\to +\infty}\sup_{t\in[0,T]}
	\left(|\uinfty'(t)-u_{n}'(t)|^{2}+
	|A^{1/2}(\uinfty(t)-u_{n}(t))|^{2}\right)=0.
	\label{th:u2uinfty}
\end{equation}

We obtain the following two results, which are probably interesting in 
themselves.

\begin{thm}[Convergence in the supercritical case]\label{thm:sup-cn2cinfty}
	Let $u_{n}(t)$ be the sequence of solutions to the approximating problems
	(\ref{pbm:eqn-cn})--(\ref{pbm:data-cn}), and let $\uinfty(t)$ be
	the solution to the limit problem
	(\ref{pbm:eqn-cinfty})--(\ref{pbm:data-cinfty}).
	
	Let us assume that
	\begin{itemize}
		\item $A$ is a self-adjoint nonnegative operator on a
		separable Hilbert space $H$,
	
		\item the coefficients $c_{n}:[0,T]\to\re$ are
		measurable and satisfy the degenerate hyperbolicity assumption
		(\ref{hp:dh}) with the same $\mu_{2}$,
		
		\item  $c_{n}(t)$ pointwise converges to $c_{\infty}(t)$, 
		
		\item $\sigma$ and $\delta$ are two positive real numbers such
		that either $\sigma>1/2$, or $\sigma=1/2$ and
		$4\delta^{2}\geq\mu_{2}$,
		
		\item $(u_{0},u_{1})\in D(A^{1/2})\times H$.
	\end{itemize}

	Then $u_{n}\to\uinfty$ in the sense of (\ref{th:u2uinfty}).
\end{thm}

\begin{thm}[Convergence in the subcritical case]\label{thm:sub-cn2cinfty}
	Let $u_{n}(t)$ be the sequence of solutions to the approximating problems
	(\ref{pbm:eqn-cn})--(\ref{pbm:data-cn}), and let $\uinfty(t)$ be
	the solution to the limit problem
	(\ref{pbm:eqn-cinfty})--(\ref{pbm:data-cinfty}).
	
	Let us assume that
	\begin{itemize}
		\item $A$ is a self-adjoint nonnegative operator on a
		separable Hilbert space $H$,
	
		\item the coefficients $c_{n}:[0,T]\to\re$ satisfy the
		strict hyperbolicity assumption (\ref{hp:sh}) with the same
		constants $\mu_{1}$ and $\mu_{2}$, and the $\omega$-continuity
		assumption (\ref{hp:ocont}) with the same continuity modulus
		$\omega(x)$,
		
		\item  $c_{n}(t)$ pointwise converges to $c_{\infty}(t)$, 
	
		\item $\sigma\in[0,1/2]$ and $\delta>0$ are two real 
		numbers satisfying (\ref{hp:sub}),
		
	\end{itemize}

	Then $u_{n}\to\uinfty$ in the sense of (\ref{th:u2uinfty}).
\end{thm}

We observe that in both cases the limit coefficient $c_{\infty}(t)$ 
satisfies the same assumptions of the approximating coefficients 
$c_{n}(t)$.

The rest of this section is devoted to the proof of these results.  A
careful inspection of the argument reveals that in both cases the
pointwise convergence assumption (\ref{hp:cn-cinfty}) can be weakened
to convergence in $L^{2}((0,T))$.

\subsection*{Estimates on components}\label{sec:conv-lemma}

In this section we provide estimates for solutions to the  
family of non-homogeneous linear ordinary differential equations
\begin{equation}
	w''(t)+2\delta\lambda^{2\sigma}w'(t)+\lambda^{2}c(t)w(t)=f(t),
	\label{ODE:f}
\end{equation}
with null initial data
\begin{equation}
	w(0)=w'(0)=0.
	\label{ODE:f-data}
\end{equation}

Our interest is motivated by the fact that Fourier components of the 
difference $u_{n}(t)-\uinfty(t)$ are solutions to problems of this 
type.

\begin{lemma}[Supercritical dissipation]\label{lemma:ode-f-sup}
	Let $T>0$, and let $w(t)$ be the solution to problem 
	(\ref{ODE:f})--(\ref{ODE:f-data}) under the following assumptions:
	\begin{itemize}
		\item the coefficient $c:[0,T]\to\re$ is measurable and
		satisfies (\ref{hp:dh}),
		
		\item $\sigma$ and $\delta$ are two positive real numbers 
		satisfying the assumptions of Theorem~\ref{thm:sup-cn2cinfty},
		
		\item $\lambda\geq 0$ and $f\in L^{2}((0,T),\re)$.
	\end{itemize}

	Then there exist two constants $\Gamma_{1}$ and $\Gamma_{2}$ such
	that
	\begin{equation}
		|w'(t)|^{2}+(1+\lambda^{2})|w(t)|^{2}\leq
		\Gamma_{1}\exp(\Gamma_{2}t)\int_{0}^{t}|f(s)|^{2}\,ds
		\quad\quad
		\forall t\in[0,T].
		\label{th:ode-f-sup}
	\end{equation}
	
	The constants $\Gamma_{1}$ and $\Gamma_{2}$ depend only on
	$\delta$, $\sigma$, and $\mu_{2}$ (in particular they are
	independent of $\lambda$, $f$, and $T$).
\end{lemma}

\paragraph{\textmd{\textit{Proof}}}

Let us consider the energy
$$E(t):=|w'(t)|^{2}+
\left(1+2\delta^{2}\lambda^{4\sigma}\right)|w(t)|^{2}
+2\delta\lambda^{2\sigma}w(t)w'(t).$$

The constants $k_{1}$, \ldots, $k_{4}$ we introduce in the sequel
are positive numbers depending only on $\delta$, $\sigma$, and
$\mu_{2}$.

Since
\begin{equation}
	|2\delta\lambda^{2\sigma}w(t)w'(t)|\leq
	\frac{3}{4}|w'(t)|^{2}+
	\frac{4}{3}\delta^{2}\lambda^{4\sigma}|w(t)|^{2},
	\label{est:lemma-misto}
\end{equation}
and since $\lambda^{2}\leq 1+\lambda^{4\sigma}$ (because $\sigma\geq
1/2$), it turns out that
\begin{equation}
	|w'(t)|^{2}+(1+\lambda^{4\sigma}+\lambda^{2})
	|w(t)|^{2}\leq k_{1}E(t).
	\label{ode-f-sup:E-equiv}
\end{equation}

The time-derivative of $E(t)$ is
\begin{eqnarray}
	E'(t) & = & -2\left(\delta\lambda^{2\sigma}|w'(t)|^{2}+
	\delta\lambda^{2\sigma+2}c(t)|w(t)|^{2}+
	\lambda^{2} c(t)w(t)w'(t)\right)
	\nonumber  \\
	\noalign{\vspace{0.5ex}}
	 &  & \mbox{}+2w'(t)f(t)+2w'(t)w(t)
	 +2\delta\lambda^{2\sigma}w(t)f(t).
	\label{ode-f-sup:E'}
\end{eqnarray}

Due to (\ref{ode-f-sup:E-equiv}), the last three terms can be estimated as
follows
\begin{eqnarray}
	\lefteqn{\hspace{-4em}
	2w'(t)f(t)+2w'(t)w(t)+2\delta\lambda^{2\sigma}w(t)f(t)} 
	\nonumber  \\
	\quad\quad & \leq & 
	2|w'(t)|^{2}+
	\left(1+\delta^{2}\lambda^{4\sigma}\right)|w(t)|^{2}+
	2|f(t)|^{2}
	\nonumber  \\
	 & \leq & k_{2}E(t)+2|f(t)|^{2}.
	\label{ode-f-sup:E'-2}
\end{eqnarray}

We claim that
\begin{equation}
	-2\left(\delta\lambda^{2\sigma}|w'(t)|^{2}+
	\delta\lambda^{2\sigma+2}c(t)|w(t)|^{2}+ \lambda^{2}
	c(t)w(t)w'(t)\right)\leq
	k_{3}E(t).
	\label{ode-f-sup:E'-1}
\end{equation}

Indeed, this inequality is equivalent to
\begin{equation}
	A|w'(t)|^{2}+
	B|w(t)|^{2}+2Cw(t)w'(t) \geq 0,
	\label{sup-quadr}
\end{equation}
where for the sake of shortness we set
$$A:=k_{3}+2\delta\lambda^{2\sigma},
\quad\quad
B:=k_{3}+2k_{3}\delta^{2}\lambda^{4\sigma}+
2\delta\lambda^{2\sigma+2}c(t),
\quad\quad
C:=k_{3}\delta\lambda^{2\sigma}+\lambda^{2}c(t).$$

The left-hand side of (\ref{sup-quadr}) is a quadratic form in $w'(t)$
and $w(t)$.  The coefficients $A$ and $B$ are positive.  Therefore,
this quadratic form is nonnegative for all values of $w'(t)$ and
$w(t)$ if and only if $AB\geq C^{2}$.  With some algebra, this
condition turns out to be equivalent to
\begin{equation}
	4\delta^{2}\lambda^{4\sigma+2}c(t)+k_{3}^{2}+
	4k_{3}\delta^{3}\lambda^{6\sigma}+
	2k_{3}\delta\lambda^{2\sigma}+
	k_{3}^{2}\delta^{2}\lambda^{4\sigma}\geq
	\lambda^{4}c^{2}(t).
	\label{est:k3}
\end{equation}

If $\sigma=1/2$ and $4\delta^{2}\geq\mu_{2}$, then
$$4\delta^{2}\lambda^{4\sigma+2}c(t)=
4\delta^{2}\lambda^{4}c(t)\geq
\mu_{2}\lambda^{4}c(t)\geq
\lambda^{4}c^{2}(t),$$
and hence even the first term in the left-hand side of (\ref{est:k3})
is greater than or equal to the right-hand side, independently of
$k_{3}$.  The same is true if $\sigma>1/2$ and 
$4\delta^{2}\lambda^{4\sigma-2}\geq\mu_{2}$.

If $\sigma>1/2$ and $4\delta^{2}\lambda^{4\sigma-2}\leq\mu_{2}$, then
$\lambda$ is bounded in terms of $\delta$ and $\mu_{2}$, and hence the
term $k_{3}^{2}$ alone in the left-hand side of (\ref{est:k3}) is
greater than or equal to the right-hand side, provided that $k_{3}$ is
suitably chosen.  This completes the proof of (\ref{ode-f-sup:E'-1}).

Plugging (\ref{ode-f-sup:E'-2}) and (\ref{ode-f-sup:E'-1}) into 
(\ref{ode-f-sup:E'}), we deduce that
$$E'(t)\leq k_{4}E(t)+2|f(t)|^{2},$$
and hence
$$E(t)\leq e^{k_{4}t}\left(E(0)+
2\int_{0}^{t}e^{-k_{4}s}|f(s)|^{2}\,ds\right)
\quad\quad
\forall t\geq 0.$$

Since $E(0)=0$ because of (\ref{ODE:f-data}), ignoring the exponential
inside the integral we conclude that
$$E(t)\leq 2e^{k_{4}t}\int_{0}^{t}|f(s)|^{2}\,ds
\quad\quad
\forall t\geq 0.$$

At this point (\ref{th:ode-f-sup}), with $\Gamma_{1}:=2k_{1}$ and 
$\Gamma_{2}:=k_{4}$, follows from (\ref{ode-f-sup:E-equiv}).\qed

\begin{lemma}[Subcritical dissipation]\label{lemma:ode-f-sub}
	Let $T>0$, and let $w(t)$ be the solution to problem 
	(\ref{ODE:f})--(\ref{ODE:f-data}) under the following assumptions:
	\begin{itemize}
		
		\item  the coefficient $c:[0,T]\to\re$ satisfies the 
		strict hyperbolicity assumption (\ref{hp:sh}) and the 
		$\omega$-continuity assumption (\ref{hp:ocont}),
	
		\item $\sigma$ and $\delta$ are two positive real numbers 
		satisfying the assumptions of Theorem~\ref{thm:sub-cn2cinfty},
			
		\item $\lambda\geq 0$ and $f\in L^{2}((0,T),\re)$.
	\end{itemize}

	Then there exist two constants $\Gamma_{3}$ and $\Gamma_{4}$ such
	that
	$$|w'(t)|^{2}+(1+\lambda^{2})|w(t)|^{2}\leq
		\Gamma_{3}\exp(\Gamma_{4}t)\int_{0}^{t}|f(s)|^{2}\,ds
		\quad\quad
		\forall t\in[0,T].$$
	
	The constants $\Gamma_{3}$ and $\Gamma_{4}$ depend only on
	$\delta$, $\sigma$, $\mu_{1}$, and on the continuity modulus
	$\omega(x)$ (in particular they are independent of $\lambda$, $f$,
	and $T$).
\end{lemma}

\paragraph{\textmd{\textit{Proof}}}

Let us extend $c(t)$ beyond $T$ by setting $c(t):=c(T)$ for $t\geq T$.
For every $\ep>0$, let us introduce the regularized coefficient
$$\cep(t):=\frac{1}{\ep}\int_{t}^{t+\ep}c(s)\,ds \quad\quad
\forall t\in[0,T].$$

It is easy to see that $\cep\in C^{1}([0,T],\re)$ and satisfies 
the following estimates:
\begin{equation}
	\mu_{1}\leq\cep(t)\leq\mu_{2}
	\quad\quad
	\forall t\in[0,T],
	\label{est:cep-sh}
\end{equation}
\begin{equation}
	|c(t)-\cep(t)|\leq\omega(\ep)
	\quad\quad
	\forall t\in[0,T],
	\label{est:cep}
\end{equation}
\begin{equation}
	|\cep'(t)|\leq\frac{\omega(\ep)}{\ep}
	\quad\quad
	\forall t\in[0,T].
	\label{est:cep'}
\end{equation}

The constants $\nu$, and $k_{1}$, \ldots, $k_{4}$ we introduce in the sequel
are positive numbers depending only on $\delta$, $\sigma$, $\mu_{1}$, and
on the continuity modulus $\omega(x)$.

To begin with, we fix $\nu\geq 1$ such that (\ref{hp:sub-lambda})
holds true.  Such a value exists due to assumption (\ref{hp:sub}). 
Then we set
$$\ep(\lambda):=\left\{
\begin{array}{ll}
	1 & \mbox{if }\lambda<\nu,  \\
	1/\lambda & \mbox{if }\lambda\geq\nu,
\end{array}
\right.$$
and we consider the energy
$$E_{\lambda}(t):=|w'(t)|^{2}+
\left(1+2\delta^{2}\lambda^{4\sigma}+
\lambda^{2}c_{\ep(\lambda)}(t)\right)|w(t)|^{2}
+2\delta\lambda^{2\sigma}w(t)w'(t).$$

The last term can be estimated as in (\ref{est:lemma-misto}). 
Therefore, keeping (\ref{est:cep-sh}) into account, it follows that
\begin{equation}
	|w'(t)|^{2}+(1+\lambda^{2})|w(t)|^{2}\leq k_{1}E_{\lambda}(t).
	\label{ode-f-sub:E-equiv}
\end{equation}

The time-derivative of $E_{\lambda}(t)$ is
\begin{eqnarray}
	E_{\lambda}'(t) & = & -2\left(\delta\lambda^{2\sigma}|w'(t)|^{2}+
	\delta\lambda^{2\sigma+2}c(t)|w(t)|^{2}+
	\lambda^{2} c(t)w(t)w'(t)\right)
	\nonumber  \\
	\noalign{\vspace{0.5ex}}
	 &  & \mbox{}+\lambda^{2}c_{\ep(\lambda)}'(t)|w(t)|^{2}
	 +2\lambda^{2}c_{\ep(\lambda)}(t)w(t)w'(t)   
	 \nonumber \\
	 \noalign{\vspace{0.5ex}}
	 & & \mbox{}+2w'(t)f(t)+2w'(t)w(t)+2\delta\lambda^{2\sigma}w(t)f(t).
	\label{ode-f-sub:E'}
\end{eqnarray}

Let $L_{1}$, $L_{2}$, $L_{3}$ denote the three lines in the 
expression of $E_{\lambda}'(t)$. To begin with, we observe that
$$L_{3}\leq 2|w'(t)|^{2}+(1+\delta^{2}\lambda^{4\sigma})|w(t)|^{2}+
2|f(t)|^{2}.$$

Since $\lambda^{4\sigma}\leq 1+\lambda^{2}$ (because $\sigma\leq
1/2$), from (\ref{ode-f-sub:E-equiv}) we deduce that
\begin{equation}
	L_{3}\leq
	k_{2}E(t)+2|f(t)|^{2}
	\quad\quad
	\forall t\in[0,T].
	\label{ode-f-sub:E'-2}
\end{equation}

Now we claim that
\begin{equation}
	L_{1}+L_{2}\leq
	k_{3}E(t)
	\quad\quad
	\forall t\in[0,T].
	\label{ode-f-sub:E'-1}
\end{equation}

Indeed, this inequality is equivalent to
\begin{equation}
	A|w'(t)|^{2}+
	B|w(t)|^{2}+2Cw(t)w'(t) \geq 0,
	\label{sub-quadr}
\end{equation}
where for the sake of shortness we set
$$A:=k_{3}+2\delta\lambda^{2\sigma},
\quad\quad\quad
C:=k_{3}\delta\lambda^{2\sigma}+
\lambda^{2}(c(t)-c_{\ep(\lambda)}(t)),$$
$$B:=2k_{3}\delta^{2}\lambda^{4\sigma}+
k_{3}+2\delta\lambda^{2\sigma+2}c(t)+
k_{3}\lambda^{2}c_{\ep(\lambda)}(t)-\lambda^{2}c_{\ep(\lambda)}'(t).$$

The left-hand side of (\ref{sub-quadr}) is a quadratic form in $w'(t)$
and $w(t)$.  The coefficient $A$ is positive.  Therefore, this
quadratic form is nonnegative for all values of $w'(t)$ and $w(t)$ if and
only if $AB\geq C^{2}$.  With some algebra, this condition turns out
to be equivalent to
\begin{equation}
	\begin{array}{c}
		\hspace{-2em}
		4\delta^{2}\lambda^{4\sigma+2}c(t)+
		4k_{3}\delta\lambda^{2\sigma+2}c_{\ep(\lambda)}(t)+k_{3}^{2} +
		4k_{3}\delta^{3}\lambda^{6\sigma}+
		2k_{3}\delta\lambda^{2\sigma}+
		k_{3}^{2}\delta^{2}\lambda^{4\sigma}+ \\
		\noalign{\vspace{1.5ex}}
		\hspace{3em}
		\mbox{}+ k_{3}^{2}\lambda^{2}c_{\ep(\lambda)}(t)\geq
		\lambda^{4}(c(t)-c_{\ep(\lambda)}(t))^{2}
		+2\delta\lambda^{2\sigma+2}c_{\ep(\lambda)}'(t) +
		k_{3}\lambda^{2}c_{\ep(\lambda)}'(t).
	\end{array}
	\label{est:sub-k3}
\end{equation}

When $\lambda\geq\nu$ we chose $\ep(\lambda)=1/\lambda$, so that 
(\ref{est:cep}) and (\ref{est:cep'}) read as
$$|c(t)-c_{\ep(\lambda)}(t)|\leq\omega\left(\frac{1}{\lambda}\right),
\hspace{3em}
|c_{\ep(\lambda)}'(t)|\leq\lambda\,\omega\left(\frac{1}{\lambda}\right).$$

Therefore, since (\ref{hp:sub-lambda}) holds true for 
$\lambda\geq\nu$, it follows that
$$\hspace{-2em}4\delta^{2}\lambda^{4\sigma+2}c(t)\geq
4\delta^{2}\mu_{1}\lambda^{4\sigma+2}\geq
\lambda^{4}\omega^{2}\left(\frac{1}{\lambda}\right)+
2\delta\lambda^{2\sigma+3}\omega\left(\frac{1}{\lambda}\right)$$
$$\hspace{3em}\geq
\lambda^{4}\left(c(t)-c_{\ep(\lambda)}(t)\right)^{2}+
2\delta\lambda^{2\sigma+2}c_{\ep(\lambda)}'(t),$$
and
$$4k_{3}\delta\lambda^{2\sigma+2}c_{\ep(\lambda)}(t)\geq
4\delta^{2}\mu_{1}\cdot\frac{k_{3}\lambda^{2\sigma+2}}{\delta}\geq
2\delta\lambda^{1-2\sigma}\omega\left(\frac{1}{\lambda}\right)
\cdot\frac{k_{3}\lambda^{2\sigma+2}}{\delta}\geq
k_{3}\lambda^{2}c_{\ep(\lambda)}'(t).$$

As a consequence, when $\lambda\geq\nu$ the sum of the first two terms
in the left-hand side of (\ref{est:sub-k3}) is greater than or equal
to the whole right-hand side, independently of $k_{3}$.

When $\lambda<\nu$ we chose $\ep(\lambda)=1$.  Thanks to
(\ref{est:cep}) and (\ref{est:cep'}) with $\ep=1$, the right-hand side
of (\ref{est:sub-k3}) is less than or equal to
$$\nu^{4}\omega^{2}(1)+2\delta\nu^{2\sigma+2}\omega(1)+
k_{3}\nu^{2}\omega(1).$$

As a consequence, the third term in the left-hand side of
(\ref{est:sub-k3}), namely $k_{3}^{2}$, is greater than or equal to
the whole right-hand side, provided that $k_{3}$ is large enough.
This completes the proof of (\ref{ode-f-sub:E'-1}).

Plugging (\ref{ode-f-sub:E'-2}) and (\ref{ode-f-sub:E'-1}) into 
(\ref{ode-f-sub:E'}), we deduce that
$$E_{\lambda}'(t)\leq k_{4}E_{\lambda}(t)+2|f(t)|^{2}
\quad\quad
\forall t\in[0,T].$$

The conclusion follows exactly as in the proof of 
Lemma~\ref{lemma:ode-f-sup}\qed

\subsection*{Proof of Theorem~\ref{thm:sup-cn2cinfty}}

The result is established in two steps.  In the first one, we prove
the conclusion under the more restrictive assumption that
$(u_{0},u_{1})\in D(A)\times D(A^{1/2})$.  In this case a stronger
result holds true, in the sense that the norm of $u_{n}-\uinfty$ in
the energy space can be estimated in terms of the norm of
$c_{n}-\cinfty$ in $L^{2}((0,T),\re)$.  In the second step we apply an
approximation procedure in order to obtain the conclusion for all
initial data, abandoning the possibility to estimate the convergence
rate.

\paragraph{\textmd{\textit{Convergence for more regular data}}}

Let us assume that $(u_{0},u_{1})\in D(A)\times D(A^{1/2})$.  We prove
that there exists a constant $\Gamma$, depending only on $\delta$, 
$\sigma$, $\mu_{2}$ and $T$, such that 
\begin{eqnarray}
	\lefteqn{\hspace{-2em}
	|\uinfty'(t)-u_{n}'(t)|^{2}+
	|A^{1/2}(\uinfty(t)-u_{n}(t))|^{2}}
	\nonumber  \\
	 & \leq & \Gamma
	\left(|A^{1/2}u_{1}|^{2}+|Au_{0}|^{2}\right)
	\int_{0}^{T}|c_{n}(s)-\cinfty(s)|^{2}\,ds
	\label{th:sup-reg-data}
\end{eqnarray}
for every $t\in[0,T]$.

To this end, we introduce the difference 
$w_{n}(t):=\uinfty(t)-u_{n}(t)$, which turns out to be a solution to
$$w_{n}''(t)+2\delta
A^{\sigma}w_{n}'(t)+\cinfty(t)Aw_{n}(t)=(c_{n}(t)-\cinfty(t))Au_{n}(t),$$
with null initial data $w_{n}(0)=w_{n}'(0)=0$. We also consider the 
components $\wnk(t)$ of $w_{n}(t)$ with respect to the orthonormal 
system $\{e_{k}\}$. If $u_{n,k}(t)$ are the corresponding components 
of $u_{n}(t)$, it turns out that $\wnk(t)$ is the solution to the 
ordinary differential equation
$$\wnk''(t)+2\delta\lk^{2\sigma}\wnk'(t)+\lk^{2}\cinfty(t)\wnk(t)=
\lk^{2}(c_{n}(t)-\cinfty(t))u_{n,k}(t),$$
with null initial data $\wnk(0)=\wnk'(0)=0$.  Therefore, we can apply
Lemma~\ref{lemma:ode-f-sup} with 
$$w(t):=\wnk(t), 
\quad\quad
\lambda:=\lk, 
\quad\quad 
c(t):=\cinfty(t), 
\quad\quad
f(t):=\lk^{2}(c_{n}(t)-\cinfty(t))u_{n,k}(t).$$

We obtain that
$$|\wnk'(t)|^{2}+\lk^{2}|\wnk(t)|^{2}\leq
\Gamma_{1}\exp(\Gamma_{2}t)\int_{0}^{t}
\lk^{4}|c_{n}(s)-\cinfty(s)|^{2}\cdot|u_{n,k}(s)|^{2}\,ds,$$
where $\Gamma_{1}$ and $\Gamma_{2}$ do not depend on $k$.  Therefore,
summing over all $k$'s we deduce that
\begin{equation}
	|w_{n}'(t)|^{2}+|A^{1/2}w_{n}(t)|^{2}\leq 
	\Gamma_{1}\exp(\Gamma_{2}t)
	\int_{0}^{t}|c_{n}(s)-\cinfty(s)|^{2}\cdot|Au_{n}(s)|^{2}\,ds.
	\label{est:wn-sup}
\end{equation}

In order to estimate $|Au_{n}(s)|$, we choose $\nu\geq 1$ satisfying
(\ref{hp:sup-nu}) (we point out that $\nu$ depends only on $\delta$,
$\sigma$ and $\mu_{2}$), and we apply Theorem~\ref{thmbibl:sup} to the
function $A^{1/2}u_{n}(t)$, which is again a solution to equation
(\ref{pbm:eqn-cn}) thanks to linearity.  We obtain that
$$|Au_{n}(t)|^{2}\leq K_{1}(\delta,\mu_{2})e^{\nu(1+\mu_{2})T}
\left(|A^{1/2}u_{1}|^{2}+|Au_{0}|^{2}\right) 
\quad\quad
\forall t\in[0,T].$$

Plugging this estimate into (\ref{est:wn-sup}), we finally deduce 
(\ref{th:sup-reg-data}).

\paragraph{\textmd{\textit{Convergence for general data}}}

Let us consider now an initial condition $(u_{0},u_{1})\in
D(A^{1/2})\times H$.  We show that, for every $\eta>0$, there exists
$n_{0}\in\n$ such that
\begin{equation}
	|u_{n}'(t)-\uinfty'(t)|^{2}+
	|A^{1/2}(u_{n}(t)-\uinfty(t))|^{2}\leq\eta
	\quad\quad
	\forall t\in[0,T],\quad\forall n\geq n_{0}.
	\label{th:sup-eta}
\end{equation}

To this end, we exploit a classical approximation argument.  We choose
$\nu\geq 1$ satisfying (\ref{hp:sup-nu}), and then we choose
$(v_{0},v_{1})\in D(A)\times D(A^{1/2})$ such that
\begin{equation}
	|u_{1}-v_{1}|^{2}+\left|A^{1/2}(u_{0}-v_{0})\right|^{2}\leq
	\frac{\eta}{9K_{1}(\delta,\mu_{2})e^{\nu(1+\mu_{2})T}},
	\label{choice-data}
\end{equation}
where $K_{1}(\delta,\mu_{2})$ is again the constant of 
Theorem~\ref{thmbibl:sup}.

Let $v_{n}(t)$ and $v_{\infty}(t)$ denote the solutions to
(\ref{pbm:eqn-cn}) and (\ref{pbm:eqn-cinfty}), respectively, with
initial data $v_{n}(0)=v_{\infty}(0)=v_{0}$ and
$v_{n}'(0)=v_{\infty}'(0)=v_{1}$.  Since $u_{n}(t)-v_{n}(t)$ is again
a solution to (\ref{pbm:eqn-cn}), from (\ref{th:sup-sumup}) and
(\ref{choice-data}) it follows that
\begin{equation}
	|u_{n}'(t)-v_{n}'(t)|^{2}+\left|A^{1/2}(u_{n}(t)-v_{n}(t))\right|^{2}
	\leq\frac{\eta}{9}
	\quad\quad
	\forall t\in[0,T].
	\label{th:sup-eta-1}
\end{equation}

Since $u_{\infty}(t)-v_{\infty}(t)$ is again a solution to
(\ref{pbm:eqn-cinfty}), from (\ref{th:sup-sumup}) and
(\ref{choice-data}) it follows that
\begin{equation}
	|u_{\infty}'(t)-v_{\infty}'(t)|^{2}+
	\left|A^{1/2}(u_{\infty}(t)-v_{\infty}(t))\right|^{2}
	\leq\frac{\eta}{9}
	\quad\quad
	\forall t\in[0,T].
	\label{th:sup-eta-2}
\end{equation}

Finally, since $(v_{0},v_{1})\in D(A)\times D(A^{1/2})$, we can apply 
(\ref{th:sup-reg-data}) to $v_{n}(t)$ and $v_{\infty}(t)$. It follows 
that there exists $n_{0}\in\n$ such that
\begin{equation}
	|v_{n}'(t)-v_{\infty}'(t)|^{2}+
	\left|A^{1/2}(v_{n}(t)-v_{\infty}(t))\right|^{2}
	\leq\frac{\eta}{9}
	\quad\quad
	\forall t\in[0,T]\quad\forall n\geq n_{0}.
	\label{th:sup-eta-3}
\end{equation}

Since
\begin{eqnarray*}
	\lefteqn{\hspace{-5em}
	|u_{n}'(t)-\uinfty'(t)|^{2}+
	\left|A^{1/2}(u_{n}(t)-\uinfty(t))\right|^{2}}  \\
	\hspace{4em} & \leq  &
	3\left(|u_{n}'(t)-v_{n}'(t)|^{2}+
	\left|A^{1/2}(u_{n}(t)-v_{n}(t))\right|^{2}\right.
	\\
	 &  & \mbox{}+|v_{n}'(t)-v_{\infty}'(t)|^{2}+
	 \left|A^{1/2}(v_{n}(t)-v_{\infty}(t))\right|^{2}  \\
	 &  & \left.\mbox{}+|u_{\infty}'(t)-v_{\infty}'(t)|^{2}+
	 \left|A^{1/2}(\uinfty(t)-v_{\infty}(t))\right|^{2}\right),
\end{eqnarray*}
conclusion (\ref{th:sup-eta}) follows from (\ref{th:sup-eta-1}), 
(\ref{th:sup-eta-2}), and (\ref{th:sup-eta-3}). This proves 
(\ref{th:u2uinfty}).\qed 

\subsection*{Proof of Theorem~\ref{thm:sub-cn2cinfty}}

Same proof of Theorem~\ref{thm:sup-cn2cinfty}, the only difference
being that the key estimates come from Theorem~\ref{thmbibl:sub} and
Lemma~\ref{lemma:ode-f-sub}, instead of Theorem~\ref{thmbibl:sup} and
Lemma~\ref{lemma:ode-f-sup}, and the constant $\Gamma$ in 
(\ref{th:sup-reg-data}) now depends on $\delta$, $\sigma$, $\mu_{1}$, 
$\mu_{2}$, $T$, and on the continuity modulus $\omega$.\qed

\setcounter{equation}{0}
\section{Interpolation spaces}\label{sec:interpolation}

In this section we introduce a family of interpolation spaces between
the energy space $D(A^{1/2})\times H$ and the standard space for
Kirchhoff equations $D(A^{3/4})\times D(A^{1/4})$.  We prove three
results.  First of all, there is propagation of regularity for
problem (\ref{pbm:lin-c(t)})--(\ref{pbm:k-data}), in the sense that
the problem is well-posed in these intermediate spaces whenever it is
well-posed in the energy space.  More important, the continuity
modulus of the function $t\to|A^{1/2}u(t)|^{2}$ depends only on the
regularity of initial data in these interpolation spaces, and not on
the regularity of $c(t)$.  Finally, we show that any pair of initial
conditions $(u_{0},u_{1})$ in the energy space lies actually in a 
suitable interpolation space.

Roughly speaking, the role of these interpolation spaces between the
energy space and $D(A^{3/4})\times D(A^{1/4})$ is the same as the role
of functions with a given continuity modulus between the space of all
continuous functions and the space of Lipschitz continuous functions.
The definition itself relies on the notion of continuity modulus.

\begin{defn}
	\begin{em}
		Let $\omega(x)$ be a continuity modulus.  For every
		$(u_{0},u_{1})\in D(A^{1/2})\times H$ we set
		$$\|(u_{0},u_{1})\|_{\omega}^{2}:=\sum_{k\in K_{0}}
		\frac{1}{\omega(1/\lk)}
		\left(|u_{1k}|^{2}+\lk^{2}|u_{0k}|^{2}\right),$$
		where $K_{0}$ denotes the set of nonnegative integers $k$ such
		that $\lk> 0$, and $u_{0k}$ and $u_{1k}$ denote the
		components of $u_{0}$ and $u_{1}$ with respect to the usual
		orthonormal system.  Then we define the space
		$$V_{\omega}:=\left\{(u_{0},u_{1})\in D(A^{1/2})\times H:
		\|(u_{0},u_{1})\|_{\omega}^{2}<+\infty\right\}.$$
	\end{em}
\end{defn}

It is easy to see that
$$D(A^{3/4})\times D(A^{1/4})\subseteq V_{\omega}\subseteq
D(A^{1/2})\times H.$$

Moreover, $V_{\omega}$ is actually a vector space, 
$\|(u_{0},u_{1})\|_{\omega}$ is a seminorm, and the full norm
$$\left(|u_{1}|^{2}+|u_{0}|^{2}+\|(u_{0},u_{1})\|_{\omega}^{2}\right)^{1/2}$$
induces a Hilbert space structure on $V_{\omega}$.

When $\omega(x)=x^{4\alpha}$ for some $\alpha\in[0,1/4]$, the space
$V_{\omega}$ is just $D(A^{\alpha+1/2})\times
D(A^{\alpha})$.

Since in this paper we are dealing with continuity moduli in several
different contexts, from now on we write $\omega_{d}(x)$ to denote
continuity moduli involved in interpolation spaces (here ``$d$\,''
stands for ``data'').

\subsection{Propagation of regularity}

When problem (\ref{pbm:lin-c(t)})--(\ref{pbm:k-data}) generates a 
continuous semigroup in the energy space, then every pair of initial 
conditions in $V_{\omega_{d}}$ gives rise to a solution lying in the 
same space for all positive times.

Let us start with the case $\sigma\geq 1/2$.

\begin{prop}[Supercritical case -- Regularity in interpolation 
	spaces]\label{prop:sup-sol-omega}
	Let $T>0$, and let $u(t)$ be the solution in $[0,T]$ to problem
	(\ref{pbm:lin-c(t)})--(\ref{pbm:k-data}) under the same assumption
	of Theorem~\ref{thmbibl:sup}.  
	
	Let us assume in addition that $(u_{0},u_{1})\in V_{\omega_{d}}$ for
	some continuity modulus $\omega_{d}(x)$.
	
	Then $(u,u')\in C^{0}([0,T],V_{\omega_{d}})$ and for every $\nu\geq 1$
	satisfying (\ref{hp:sup-nu}) it turns out that
	\begin{equation}
		\|(u(t),u'(t))\|_{\omega_{d}}^{2}\leq 
		K_{1}(\delta,\mu_{2})e^{\nu(1+\mu_{2}) t}
		\|(u_{0},u_{1})\|_{\omega_{d}}^{2}
		\quad\quad
		\forall t\in[0,T],
		\label{th:od-cont}
	\end{equation}
	where $K_{1}(\delta,\mu_{2})$ is the constant of 
	Theorem~\ref{thmbibl:sup}.
\end{prop}

\paragraph{\textmd{\textit{Proof}}}

Let $\{u_{k}(t)\}$ denote the components of $u(t)$ with respect to the
usual orthonormal system $\{e_{k}\}$, and let $\{u_{0k}\}$ and
$\{u_{1k}\}$ denote the corresponding components of initial
conditions. Since we are in the assumptions of 
Theorem~\ref{thmbibl:sup}, we can estimate these components as in 
(\ref{est:uk-}) and (\ref{est:uk+sup}). Since 
$K_{1}(\delta,\mu_{2})\geq 1$, we can  combine these two estimates 
and deduce that
$$|u_{k}'(t)|^{2}+\lk^{2}|u_{k}(t)|^{2}\leq
K_{1}(\delta,\mu_{2})e^{\nu(1+\mu_{2})t}
\left(|u_{1k}|^{2}+\lk^{2}|u_{0k}|^{2}\right)
\quad\quad
\forall t\in[0,T],$$
independently of $k$.  Dividing by $\omega(1/\lk)$, and summing over
$K_{0}$ (the set of indices $k$ with $\lk\neq 0$), we obtain
(\ref{th:od-cont}).

The same estimate shows also the uniform convergence of the series
$$\sum_{k\in K_{0}}\frac{\lk^{2}}{\omega_{d}(1/\lk)}|u_{k}(t)|^{2},
\hspace{4em}
\sum_{k\in K_{0}}\frac{1}{\omega_{d}(1/\lk)}|u_{k}'(t)|^{2},$$
which proves the continuity of the pair $(u,u')$ with values in
$V_{\omega_{d}}$.\qed
\medskip

The result for the subcritical case is analogous. Here we state 
explicitly a time independent estimate for high-frequency components 
in the same spirit of (\ref{th:sub-u+}). This estimate is crucial in 
the proof of the global existence statement of 
Theorem~\ref{thm:sub-K}. An analogous estimate holds true also in the 
supercritical case, but we do not need it in the sequel.

\begin{prop}[Subcritical case -- Regularity in interpolation spaces]
	\label{prop:sub-sol-omega}
	Let $T>0$, and let $u(t)$ be the solution in $[0,T]$ to problem
	(\ref{pbm:lin-c(t)})--(\ref{pbm:k-data}) under the same assumption
	of Theorem~\ref{thmbibl:sub}.  
	
	Let us assume in addition that $(u_{0},u_{1})\in V_{\omega_{d}}$ for
	some continuity modulus $\omega_{d}(x)$.
	
	Then $(u,u')\in C^{0}([0,T],V_{\omega_{d}})$ and for every $\nu\geq 1$
	satisfying (\ref{hp:sub-lambda}) it turns out that
	\begin{equation}
		\|(u(t),u'(t))\|_{\omega_{d}}^{2}\leq 
		K_{2}(\delta,\mu_{1},\mu_{2})e^{\nu(1+\mu_{2}) t}
		\|(u_{0},u_{1})\|_{\omega_{d}}^{2}
		\quad\quad
		\forall t\in[0,T],
		\label{th:od-sub-cont}
	\end{equation}
	where $K_{2}(\delta,\mu_{1},\mu_{2})$ is the constant of 
	Theorem~\ref{thmbibl:sub}.
	
	Moreover, the high-frequency component of $u(t)$ satisfies
	\begin{equation}
		\|(u_{\nu,+}(t),u_{\nu,+}'(t))\|_{\omega_{d}}^{2}\leq 
		K_{2}(\delta,\mu_{1},\mu_{2})
		\|(u_{0,\nu,+},u_{1,\nu,+})\|_{\omega_{d}}^{2}
		\quad\quad
		\forall t\in[0,T].
		\label{th:od-sub-cont+}
	\end{equation}
\end{prop}

\paragraph{\textmd{\textit{Proof}}}

Same proof of Proposition~\ref{prop:sup-sol-omega}, this time starting
from estimates (\ref{est:uk-}) and (\ref{est:uk+sub}), which hold true
under the assumptions of Theorem~\ref{thmbibl:sub}.\qed

\begin{rmk}
	\begin{em}
		From (\ref{hp:omega-monot-0}) it follows that
		$\omega_{d}(1/\lk)\geq\omega_{d}(1/\nu)$ whenever
		$\lk\leq\nu$.  Therefore, both in the supercritical and in the
		subcritical case, the low-frequency components can be
		estimated in terms of the usual energy as follows
		\begin{equation}
			\|(u_{\nu,-}(t),u_{\nu,-}'(t))\|_{\omega_{d}}^{2}\leq 
			\frac{1}{\omega_{d}(1/\nu)}\left(
			|u_{\nu,-}'(t)|^{2}+|A^{1/2}u_{\nu,-}(t)|^{2}\right)
			\quad\quad
			\forall t\in[0,T].
			\label{th:od-sub-cont-}
		\end{equation}
		Also this estimate is crucial in the proof of the global
		existence statement of Theorem~\ref{thm:sub-K}.  
	\end{em}
\end{rmk}

\subsection{Interpolation spaces vs time regularity}

It is well-known that space regularity and time regularity of
solutions to hyperbolic problems are strongly related.  The following
result clarifies the connection between the interpolation space
$V_{\omega_{d}}$ and the continuity modulus of the function
$t\to|A^{1/2}u(t)|^{2}$.

\begin{prop}[Continuity modulus of $|A^{1/2}u(t)|^{2}$]\label{prop:o-cont-auq}
	Let $\omega_{d}(x)$ be a continuity modulus, let $T>0$, and let
	$u:[0,T]\to H$ be a function with the regularity 
	(\ref{th:u-reg}).
	
	Let us assume that $(u(t),u'(t))\in V_{\omega_{d}}$ for every 
	$t\in[0,T]$, and there exists a constant $L$ such that
	\begin{equation}
		\|(u(t),u'(t))\|_{\omega_{d}}^{2}\leq L
		\quad\quad
		\forall t\in[0,T].
		\label{hp:g-omega-cont}
	\end{equation}
	
	Then it turns out that
	\begin{equation}
		\left||A^{1/2}u(a)|^{2}-|A^{1/2}u(b)|^{2}\right|\leq
		3L\,\omega_{d}(|a-b|)
		\quad\quad
		\forall (a,b)\in[0,T]^{2}.
		\label{th:g-omega-cont}
	\end{equation}
\end{prop}

\paragraph{\textmd{\textit{Proof}}}

Let $u_{k}(t)$ denote the components of $u(t)$ with respect to the
usual orthonormal system.  For every $\ep>0$, let $K_{\ep}^{-}$ denote
the set of nonnegative integers $k$ such that $0<\lk\leq 1/\ep$, and
let $K_{\ep}^{+}$ denote the set of nonnegative integers $k$ such that
$\lk>1/\ep$.  Let us set for simplicity
$$g(t):=|A^{1/2}u(t)|^{2}=\sum_{k=0}^{\infty}\lk^{2}|u_{k}(t)|^{2},
\hspace{4em}
\gep(t):=\sum_{k\in K_{\ep}^{-}}\lk^{2}|u_{k}(t)|^{2}.$$

The time-derivative of $\gep(t)$ exists and is given by
$$\gep'(t)=2\sum_{k\in K_{\ep}^{-}}\lk^{2}u_{k}(t)u_{k}'(t),$$
and hence
\begin{equation}
	|\gep'(t)|\leq\sum_{k\in K_{\ep}^{-}}\lk
	\left(|u_{k}'(t)|^{2}+\lk^{2}|u_{k}(t)|^{2}\right).
	\label{est:gep'}
\end{equation}

For every $k\in K_{\ep}^{-}$ it turns out that $1/\lk\geq\ep$, and
hence $\lk\,\omega_{d}(1/\lk)\leq\omega_{d}(\ep)/\ep$ because of the
monotonicity property (\ref{hp:omega-monot}) of the continuity
modulus.  Therefore, from (\ref{est:gep'}) and assumption 
(\ref{hp:g-omega-cont}) we deduce that
$$|\gep'(t)|\leq\frac{\omega_{d}(\ep)}{\ep}\sum_{k\in K_{\ep}^{-}}
\frac{1}{\omega_{d}(1/\lk)}
\left(|u_{k}'(t)|^{2}+\lk^{2}|u_{k}(t)|^{2}\right)\leq
L\,\frac{\omega_{d}(\ep)}{\ep}$$
for every $t\in[0,T]$.  In particular, from the mean value theorem it
follows that
\begin{equation}
	|\gep(a)-\gep(b)|\leq |a-b|\cdot\max_{t\in[0,T]}|\gep'(t)|
	\leq L\,\frac{\omega_{d}(\ep)}{\ep}\cdot|a-b|
	\label{est:gep-diff}
\end{equation}
for every $a$ and $b$ in $[0,T]$.

On the contrary, for every $k\in K_{\ep}^{+}$ it turns out that
$1/\lk<\ep$, and hence $\omega_{d}(1/\lk)\leq\omega_{d}(\ep)$ because of the
monotonicity property (\ref{hp:omega-monot-0}) of the continuity
modulus.  It follows that
$$\hspace{-3em}|g(t)-\gep(t)|=\sum_{k\in
K_{\ep}^{+}}\lk^{2}|u_{k}(t)|^{2} \leq\sum_{k\in K_{\ep}^{+}}
\left(|u_{k}'(t)|^{2}+\lk^{2}|u_{k}(t)|^{2}\right)$$
$$\hspace{3em}\leq
\omega_{d}(\ep)\sum_{k\in K_{\ep}^{+}}\frac{1}{\omega_{d}(1/\lk)}
\left(|u_{k}'(t)|^{2}+\lk^{2}|u_{k}(t)|^{2}\right)
\leq L\,\omega_{d}(\ep)$$
for every $t\in[0,T]$.  From this estimate and (\ref{est:gep-diff}) we
deduce that
\begin{eqnarray*}
	|g(a)-g(b)| & \leq & |g(a)-g_{\ep}(a)|+
	|g_{\ep}(a)-g_{\ep}(b)|+|g_{\ep}(b)-g(b)|\\
	\noalign{\vspace{1ex}}
	 & \leq & L\,\omega_{d}(\ep)+
	 L\,\frac{\omega_{d}(\ep)}{\ep}\cdot|a-b|+
	 L\,\omega_{d}(\ep)
\end{eqnarray*}
for every $a$ and $b$ in $[0,T]$, and every $\ep>0$. Now it
is enough to choose $\ep:=|a-b|$, and (\ref{th:g-omega-cont}) is
proved.\qed

\subsection{An intermediate space for each initial condition}

In the next result we show that every pair $(u_{0},u_{1})$ which is in
$D(A^{\alpha+1/2})\times D(A^{\alpha})$ for some $\alpha\in[0,1/4)$
lies actually in a ``better'' space, namely a space $V_{\omega_{d}}$
corresponding to a continuity modulus $\omega_{d}(x)$ which tends to 0
faster than $x^{4\alpha}$.  It is a refinement of the classical calculus
result according to which the terms of a converging series can always
be multiplied by a diverging sequence obtaining again a converging
series. We can also control the (semi)norm in the interpolation space in 
terms of the (semi)norm in the original space.

\begin{prop}[Customized interpolation spaces]\label{prop:data-omega}
	For every $\alpha\in[0,1/4)$ and every $(u_{0},u_{1})\in 
	D(A^{\alpha+1/2})\times D(A^{\alpha})$, there exists a continuity 
	modulus $\omega_{d}(x)$ such that $(u_{0},u_{1})\in V_{\omega_{d}}$ and
	\begin{equation}
		\omega_{d}(1)=1,
		\label{th:od=1}
	\end{equation}
	\begin{equation}
		\omega_{d}(x)\leq x^{4\alpha}
		\quad\quad
		\forall x\geq 0,
		\label{th:od<=}
	\end{equation}
	\begin{equation}
		\lim_{x\to 0^{+}}\frac{\omega_{d}(x)}{x^{4\alpha}}=0,
		\label{th:od-lim}
	\end{equation}
	\begin{equation}
		\|(u_{0},u_{1})\|_{\omega_{d}}^{2}\leq 2\left(
		|A^{\alpha}u_{1}|^{2}+|A^{\alpha+1/2}u_{0}|^{2}\right).
		\label{th:od-data}
	\end{equation}
\end{prop}

\paragraph{\textmd{\textit{Proof}}}

Let $u_{0k}$ and $u_{1k}$ denote the components of $u_{0}$ and $u_{1}$
with respect to the usual orthonormal system.  Let us set for
simplicity
$$E_{k}:=\lk^{4\alpha}\left(|u_{1k}|^{2}+\lk^{2}|u_{0k}|^{2}\right)$$
and
$$E:=|A^{\alpha}u_{1}|^{2}+|A^{\alpha+1/2}u_{0}|^{2}=
\sum_{k=0}^{\infty}E_{k}.$$

For every $n\in\n$, let us set
$$A_{n}:=\left\{k\in\n:n\leq\lk<n+1\right\},
\hspace{3em}
a_{n}:=\sum_{k\in A_{n}}E_{k}.$$

Since 
$$\sum_{n=0}^{\infty}a_{n}=\sum_{k=0}^{\infty}E_{k}=E,$$
there exists an increasing sequence $n_{h}$ of nonnegative integers 
such that $n_{0}=0$ and
\begin{equation}
	\sum_{n\geq n_{h}}a_{n}\leq\frac{E}{4^{h}}
	\quad\quad
	\forall h\in\n.
	\label{inutile}
\end{equation}

For every $h\in\n$, let us set
$$B_{h}:=\left\{k\in\n:n_{h}\leq\lk<n_{h+1}\right\}.$$

Since $B_{h}\subseteq\cup_{n\geq n_{h}}A_{n}$, from (\ref{inutile}) it
turns out that
\begin{equation}
	\sum_{k\in B_{h}}E_{k}\leq\sum_{n\geq n_{h}}a_{n}
	\leq\frac{E}{4^{h}}
	\quad\quad
	\forall h\in\n.
	\label{est:sum-Bh}
\end{equation}

Let us consider the sequence $\varphi_{h}$ defined by 
$\varphi_{0}=1$,  $\varphi_{1}=1$, 
and
$$\varphi_{h+1}=\min\left\{2^{h},
\frac{n_{h+1}}{n_{h}}\varphi_{h}\right\}
\quad\quad
\forall h\geq 1.$$

From this definition it follows that
\begin{equation}
	\varphi_{h+1}\leq 2^{h},
	\hspace{3em}
	\varphi_{h+1}\geq\varphi_{h},
	\hspace{3em}
	\frac{n_{h+1}}{\varphi_{h+1}}\geq\frac{n_{h}}{\varphi_{h}}
	\label{th:phik}
\end{equation}
for every $h\in\n$. Moreover it turns out that
\begin{equation}
	\lim_{h\to +\infty}\varphi_{h}=+\infty.
	\label{th:phik-lim}
\end{equation}

Indeed, (\ref{th:phik-lim}) is obvious if $\varphi_{h+1}=2^{h}$ for 
infinitely many indices. If not, it means that 
$\varphi_{h+1}=\varphi_{h}n_{h+1}/n_{h}$ for every $h$ greater than 
or equal to some $h_{0}\geq 1$. In this case an easy induction gives 
that
$$\varphi_{h}=\frac{n_{h}}{n_{h_{0}}}\varphi_{h_{0}}
\quad\quad
\forall h\geq h_{0},$$
so that (\ref{th:phik-lim}) follows from the fact that $n_{h}\to 
+\infty$.

Let us consider now the piecewise affine function
$\varphi:[0,+\infty)\to\re$ such that $\varphi(n_{h})=\varphi_{h}$ for
every $h\in\n$, namely the function defined by
$$\varphi(x):=\varphi_{h}+
\dfrac{\varphi_{h+1}-\varphi_{h}}{n_{h+1}-n_{h}}(x-n_{h})
\quad\quad 
\mbox{if }n_{h}\leq x\leq n_{h+1}\mbox{ for some }h\in\n.$$

From the first relation in (\ref{th:phik}) it follows that
\begin{equation}
	\varphi(x)\leq 2^{h}
	\quad\quad
	\forall x\in[n_{h},n_{h+1}].
	\label{th:phi-2h}
\end{equation}

Moreover $\varphi(1)=1$ and
\begin{equation}
	\varphi(x)\geq 1
	\quad\quad
	\forall x\geq 0.
	\label{th:phi>=1}
\end{equation}

From the second relation in (\ref{th:phik}) it follows that 
$\varphi(x)$ is nondecreasing, and from (\ref{th:phik-lim}) it 
follows that
\begin{equation}
	\lim_{x\to +\infty}\varphi(x)=+\infty.
	\label{th:phi-lim}
\end{equation}

Finally, from the third relation in (\ref{th:phik}) one can prove 
that
\begin{equation}
	\mbox{the function }x\to\frac{x}{\varphi(x)}\mbox{ is 
	nondecreasing}
	\label{th:phi-monot}
\end{equation}
(it is enough to show that its derivative is nonnegative in each 
interval $(n_{h},n_{h+1})$).

Let us finally set
$$\omega_{d}(x):=\left\{
\begin{array}{ll}
	0 & \mbox{if }x=0,  \\
	\dfrac{x^{4\alpha}}{[\varphi(1/x)]^{1-4\alpha}}\quad & \mbox{if 
	}x>0.
\end{array}
\right.$$

We claim that $\omega_{d}(x)$ is a continuity modulus, and that
(\ref{th:od=1}) through (\ref{th:od-data}) hold true.  Equality
(\ref{th:od=1}) follows from the fact that $\varphi(1)=1$, while
estimate (\ref{th:od<=}) follows from (\ref{th:phi>=1}).  Since
$4\alpha<1$, from (\ref{th:phi-lim}) we deduce that $\omega_{d}(x)\to
0$ as $x\to 0^{+}$, which proves that $\omega_{d}$ is continuous also
in $x=0$, the only point in which continuity was nontrivial.  The
limit (\ref{th:od-lim}) follows from (\ref{th:phi-lim}) because
$4\alpha<1$.  The monotonicity property (\ref{hp:omega-monot-0}) of
$\omega_{d}(x)$ follows from the fact that $\varphi(x)$ is
nondecreasing, while the monotonicity property (\ref{hp:omega-monot})
is equivalent to (\ref{th:phi-monot}) after the variable change
$y:=1/x$.

It remains to prove (\ref{th:od-data}). To this end, we begin by 
observing that
\begin{eqnarray*}
	\|(u_{0},u_{1})\|_{\omega_{d}}^{2} & = &
	\sum_{k\in K_{0}}\frac{1}{\omega_{d}(1/\lk)}
	\left(|u_{1k}|^{2}+\lk^{2}|u_{0k}|^{2}\right) \\
	 & \leq &
	\sum_{k=0}^{\infty}\left[\varphi(\lk)\right]^{1-4\alpha}E_{k} \\
	 & = & \sum_{h=0}^{\infty}\sum_{k\in B_{h}}
	\left[\varphi(\lk)\right]^{1-4\alpha}E_{k}.  
\end{eqnarray*}

Now for every $k\in B_{h}$ it turns out that $n_{h}\leq\lk<n_{h+1}$. 
From (\ref{th:phi-2h}) and the fact that $4\alpha<1$ it follows that
$$\left[\varphi(\lk)\right]^{1-4\alpha}\leq
2^{h(1-4\alpha)}\leq 2^{h}.$$

Keeping (\ref{est:sum-Bh}) into account, we deduce that
$$\|(u_{0},u_{1})\|_{\omega_{d}}^{2}\leq
\sum_{h=0}^{\infty}\sum_{k\in B_{h}}
\left[\varphi(\lk)\right]^{1-4\alpha}E_{k}\leq
\sum_{h=0}^{\infty}2^{h}\sum_{k\in B_{h}}E_{k}\leq
\sum_{h=0}^{\infty}2^{h}\frac{E}{4^{h}}=2E,$$
which proves (\ref{th:od-data}).\qed

\setcounter{equation}{0}
\section{A priori estimates for Kirchhoff equations}\label{sec:a-priori}

In this section we derive a priori estimates for solutions to 
(\ref{pbm:k-eqn}). Classical estimates are based on the usual 
Hamiltonian
\begin{equation}
	H(t):=|u'(t)|^{2}+M\left(|A^{1/2}u(t)|^{2}\right),
	\label{defn:H}
\end{equation}
where
$$M(x):=\int_{0}^{x}m(s)\,ds
\quad\quad
\forall x\geq 0.$$

If $u(t)$ is a solution to (\ref{pbm:k-eqn})--(\ref{pbm:k-data}) in 
some time-interval $[0,T]$, then the formal time-derivative of $H(t)$ 
is given by
\begin{equation}
	H'(t)=-4\delta\left|A^{\sigma/2}u'(t)\right|^{2},
	\label{deriv:H'}
\end{equation}
and therefore
\begin{equation}
	|u'(t)|^{2}+M\left(|A^{1/2}u(t)|^{2}\right)+
	4\delta\int_{0}^{t}\left|A^{\sigma/2}u'(s)\right|^{2}\,ds=
	|u_{1}|^{2}+M\left(|A^{1/2}u_{0}|^{2}\right)
	\label{est:ap-easy}
\end{equation}
for every $t\in[0,T]$.  This is the classical way to obtain estimates
both on $|u'(t)|$ and on $|A^{1/2}u(t)|$.  Unfortunately, things are
not so simple under our assumptions.

A first problem comes from the lack of strict hyperbolicity when
$\sigma>1/2$.  In that case there is no guarantee that $M(x)\to
+\infty$ as $x\to +\infty$, and hence (\ref{est:ap-easy}) does not
provide a bound on $|A^{1/2}u(t)|$, not even locally.

The second and even worse problem comes from regularity issues.
Indeed, when computing the time-derivative of $H(t)$, we have to deal
with terms such as $\langle u'(t),Au(t)\rangle$, which are quite 
delicate in the case of weak solutions. In the literature this issue 
has been addressed in two different ways.
\begin{itemize}
	\item If we write the term in the form $\langle
	A^{1/4}u'(t),A^{3/4}u(t)\rangle$, then it makes sense for
	solutions living in the phase space $D(A^{3/4})\times D(A^{1/4})$,
	which is indeed the classical space both for Kirchhoff equations
	without dissipation, and for Kirchhoff equations with standard
	dissipation ($\sigma=0$).  This argument is probably hopeless in
	the case of less regular solutions.  For example, when $\delta=0$
	equation (\ref{deriv:H'}) seems to suggest that $H(t)$ is constant
	along trajectories, but as far as we know there is no rigorous
	proof of this fact for solutions with regularity less than
	$D(A^{3/4})\times D(A^{1/4})$.

	\item In the strongly dissipative case with $\sigma=1$, equality
	(\ref{est:ap-easy}) seems to suggest that $A^{1/2}u'(t)$ lies in
	$L^{2}((0,T),H)$.  Therefore, if we write the term $\langle
	u'(t),Au(t)\rangle$ in the form $\langle
	A^{1/2}u'(t),A^{1/2}u(t)\rangle$, then it makes sense (at least
	almost everywhere) also for solutions in the energy space.  This
	is the key point in the paper~\cite{matos-pereira}, in which
	Nishihara's theory~\cite{nishihara:k-strong} is extended to
	initial data in the energy space.  Unfortunately, this heuristic
	argument strongly suggests that this approach is hopeless when
	$\sigma<1$.
\end{itemize}

Our assumptions are apparently too weak for both strategies.  This
leads us to follow a different path.  
As for the lack of coerciveness of $M(x)$, we introduce a modified
Hamiltonian.
As for the lack of regularity of solutions, when $\sigma\geq 1$ we
exploit a strategy similar to the one described in the second point
above.  When $\sigma\in(0,1)$, we exploit the regularizing effect
presented in section~\ref{sec:regularizing}.  

In the first a priori estimate we control $|A^{1/2}u(t)|$ in the
supercritical case.  We obtain an exponential bound, which is far from
being optimal and could be improved to a polynomial bound with some
additional effort.  Nevertheless, any bound is enough in the sequel.

\begin{prop}[Supercritical dissipation -- A priori 
	bound]\label{prop:ap-sup}
	Let us consider problem (\ref{pbm:k-eqn})--(\ref{pbm:k-data})
	under the same assumptions of Theorem~\ref{thm:sup-K}.  
	
	Let us assume that there exists a local solution $u(t)$, defined
	in the time-interval $[0,T]$ for some $T>0$, with the regularity
	prescribed by (\ref{th:u-reg}).
	
	Then there exists two constants $K_{3}(\delta)$ and 
	$K_{4}(\delta)$, depending only on $\delta$ (an in particular 
	independent of $T$ and $u$), such that
	\begin{equation}
		|A^{1/2}u(t)|^{2}\leq K_{3}(\delta)\left(
		|u_{1}|^{2}+|u_{0}|^{2}+|A^{1/2}u_{0}|^{2}
		+M\left(|A^{1/2}u_{0}|^{2}\right)\right)
		\exp(K_{4}(\delta)t)
		\label{th:ap-sup-k}
	\end{equation}
	for every $t\in[0,T]$.
\end{prop}

\paragraph{\textmd{\textit{Proof}}}

Let us distinguish the cases $\sigma\geq 1$ and $1/2<\sigma<1$.

\subparagraph{\textmd{\textit{Case $\sigma\geq 1$}}}

Let us set
\begin{equation}
	c(t):=m\left(|A^{1/2}u(t)|^{2}\right).
	\label{defn:c(t)}
\end{equation}

Let $\{u_{k}(t)\}$ denote the components of $u(t)$ with respect to the 
usual orthonormal system $\{e_{k}\}$, and let
$$v_{n}(t):=\sum_{k=0}^{n}u_{k}(t)e_{k},
\quad\quad\quad
c_{n}(t):=m\left(|A^{1/2}v_{n}(t)|^{2}\right).$$

It is easy to see that $v_{n}(t)$ is a solution to the linear problem
$$v_{n}''(t)+2\delta A^{\sigma}v_{n}'(t)+c(t)Av_{n}(t)=0.$$

From the regularity assumption (\ref{th:u-reg}) it turns out that
$$v_{n}'(t)\to u'(t),
\quad\quad
A^{1/2}v_{n}(t)\to A^{1/2}u(t),
\quad\quad
c_{n}(t)\to c(t),$$
and the convergence is uniform in $[0,T]$ in all three cases (but in
the sequel we need just pointwise convergence and uniform
boundedness).  Let us consider the modified approximated Hamiltonian
$$H_{n}(t):=|v_{n}'(t)|^{2}+|A^{1/2}v_{n}(t)|^{2}+
M\left(|A^{1/2}v_{n}(t)|^{2}\right).$$

Its time-derivative is
$$H_{n}'(t)=-4\delta|A^{\sigma/2}v_{n}'(t)|^{2}+
2(c_{n}(t)-c(t)+1)\langle A^{1/2}v_{n}(t), A^{1/2}v_{n}'(t)\rangle.$$

Since $\sigma\geq 1$, it turns out that
$$|A^{1/2}v_{n}'(t)|^{2}\leq
|A^{\sigma/2}v_{n}'(t)|^{2}+|v_{n}'(t)|^{2},$$
and hence
\begin{eqnarray*}
	H_{n}'(t) & \leq & -4\delta|A^{\sigma/2}v_{n}'(t)|^{2}+
	\frac{1}{4\delta}(|c_{n}(t)-c(t)|+1)^{2}|A^{1/2}v_{n}(t)|^{2}+
	4\delta|A^{1/2}v_{n}'(t)|^{2}\\
	 & \leq & \frac{1}{4\delta}(|c_{n}(t)-c(t)|+1)^{2}|A^{1/2}v_{n}(t)|^{2}+
	4\delta|v_{n}'(t)|^{2}  \\
	 & \leq & \left(\frac{1}{4\delta}(|c_{n}(t)-c(t)|+1)^{2}
	 +4\delta\right)H_{n}(t).
\end{eqnarray*}

Integrating this differential inequality we obtain that
$$|A^{1/2}v_{n}(t)|^{2}\leq H_{n}(t)\leq H_{n}(0)
\exp\left(\frac{1}{4\delta}\int_{0}^{t}
(|c_{n}(s)-c(s)|+1)^{2}\,ds+4\delta t\right).$$

Passing to the limit as $n\to +\infty$ we conclude that
$$|A^{1/2}u(t)|^{2}\leq 
\left(|u_{1}|^{2}+|A^{1/2}u_{0}|^{2}+M\left(|A^{1/2}u_{0}|^{2}\right)\right)
\exp\left(\frac{t}{4\delta}+4\delta t\right),$$
which proves (\ref{th:ap-sup-k}) in the case $\sigma\geq 1$ with
$K_{3}(\delta):=1$ and
$K_{4}(\delta):=4\delta+(4\delta)^{-1}$.

\subparagraph{\textmd{\textit{Case $1/2<\sigma<1$}}}

Let us consider $u(t)$ as a solution to the linear equation
(\ref{pbm:lin-c(t)}) with the coefficient $c(t)$ defined as in
(\ref{defn:c(t)}).  Let us consider the Hamiltonian $H(t)$ defined 
in (\ref{defn:H}). Since $\sigma\in(1/2,1)$, from 
Theorem~\ref{thmbibl:sup-gevrey} we know that 
$u'(t)$ is continuous with values in $D(A^{1/2})$ for every 
$t\in(0,T]$, or at least for every $t\in(0,S_{0})\cup(S_{0},T]$ for 
some $S_{0}\in(0,T)$. 

This is enough to justify rigorously the calculation leading to
(\ref{deriv:H'}).  Indeed, the terms such as $\langle
u'(t),Au(t)\rangle$ which appear during the computation of $H'(t)$ can
be interpreted in the form $\langle A^{1/2}u'(t),A^{1/2}u(t)\rangle$,
and hence they are well-defined when both $u(t)$ and $u'(t)$ are
continuous with values in $D(A^{1/2})$.  This proves that
(\ref{deriv:H'}) holds true for every $t\in(0,T]$, or at least for
every $t\in(0,S_{0})\cup(S_{0},T]$.  Since $H(t)$ is continuous
because of the regularity assumption (\ref{th:u-reg}), this is enough
to conclude that (\ref{est:ap-easy}) holds true for every $t\in[0,T]$.

Let us consider now the modified Hamiltonian
$$\Hhat(t)=H(t)+\delta^{2}|u(t)|^{2}+
\delta^{2}|A^{1/2}u(t)|^{2}+
\delta\langle A^{1-\sigma}u(t),u'(t)\rangle,$$
which is well-defined in the range $1/2\leq\sigma\leq 1$.  In this 
range it turns out that
\begin{eqnarray*}
	\left|\delta\langle A^{1-\sigma}u(t),u'(t)\rangle\right| & \leq & 
	\delta|A^{1-\sigma}u(t)|\cdot|u'(t)|  \\
	 & \leq & \frac{1}{2}|u'(t)|^{2}+\frac{\delta^{2}}{2}|A^{1-\sigma}u(t)|^{2}  \\
	 & \leq & \frac{1}{2}|u'(t)|^{2}+\frac{\delta^{2}}{2}|A^{1/2}u(t)|^{2}+
	 \frac{\delta^{2}}{2}|u(t)|^{2},
\end{eqnarray*}
and therefore
\begin{equation}
	\Hhat(0)\leq 2|u_{1}|^{2}+
	2\delta^{2}|u_{0}|^{2}+
	2\delta^{2}|A^{1/2}u_{0}|^{2}+
	M\left(|A^{1/2}u_{0}|^{2}\right)
	\label{est:Hhat(0)}
\end{equation}
and
\begin{equation}
	\frac{\delta^{2}}{2}|A^{1/2}u(t)|^{2}\leq\Hhat(t)
	\quad\quad
	\forall t\in[0,T].
	\label{est:a1/2-hhat}
\end{equation}

Now we claim that $\Hhat(t)$ is time-differentiable where $H(t)$ is
time-differentiable, namely in $(0,T]$ or in $(0,S_{0})\cup(S_{0},T]$,
and its time-derivative is
$$\Hhat'(t) =
-4\delta|A^{\sigma/2}u'(t)|^{2}+ 2\delta^{2}\langle u(t),u'(t)\rangle
+ \delta|A^{(1-\sigma)/2}u'(t)|^{2}- 
\delta c(t)|A^{1-\sigma/2}u(t)|^{2}.$$

Indeed, in addition to the terms due to $H'(t)$, now the most
dangerous term is the last one.  In order to show that this term is
well-defined, we consider the same three cases as in the statement of
Theorem~\ref{thmbibl:sup-gevrey}.  In the first case, the term is
well-defined because of the regularity of $u(t)$.  In the second case,
it is well-defined because $c(t)\equiv 0$ in $[0,T]$.  In the third
case, it is well-defined in $[0,S_{0}]$ because $c(t)\equiv 0$ in that
interval, and it is well-defined in $(S_{0},T]$ because of the
regularity of $u(t)$ in that interval.

In order to estimate $\Hhat'(t)$, we observe that
$$2\delta^{2}\langle u(t),u'(t)\rangle\leq
\delta^{2}|u(t)|^{2}+\delta^{2}|u'(t)|^{2},$$
and, since $\sigma\geq 1/2$,
$$\delta|A^{(1-\sigma)/2}u'(t)|^{2}\leq
\delta|A^{\sigma/2}u'(t)|^{2}+\delta|u'(t)|^{2}.$$

It follows that
$$\Hhat'(t)\leq\delta^{2}|u(t)|^{2}+(\delta^{2}+\delta)|u'(t)|^{2}.$$

Now we estimate $|u'(t)|$ by means of (\ref{est:ap-easy}), and we
estimate $|u(t)|$ by observing that 
$$|u(t)|\leq
|u_{0}|+\int_{0}^{t}|u'(s)|\,ds\leq
|u_{0}|+tH(0)^{1/2},$$
and hence
$$|u(t)|^{2}\leq 2|u_{0}|^{2}+2H(0)t^{2}.$$

It follows that
$$\Hhat'(t)\leq 2\delta^{2}|u_{0}|^{2}+2\delta^{2}H(0)t^{2}+
(\delta^{2}+\delta)H(0)$$
for every $t\in(0,T]$, or at least for every
$t\in(0,S_{0})\cup(S_{0},T]$. Since $\Hhat(t)$ is continuous in 
$[0,T]$, a simple integration gives that
$$\Hhat(t)\leq\Hhat(0)+2\delta^{2}|u_{0}|^{2}t+
(\delta^{2}+\delta)H(0)t+
2\delta^{2}H(0)\frac{t^{3}}{3}
\quad\quad
\forall t\in[0,T].$$

Now we estimate $\Hhat(0)$ as in (\ref{est:Hhat(0)}) and
powers of $t$ with $e^{t}$.  Keeping (\ref{est:a1/2-hhat})
into account, we obtain an estimate of the form (\ref{th:ap-sup-k}), 
with $K_{4}(\delta):=1$ and a suitable $K_{3}(\delta)$. 

Finally, if we want $K_{3}(\delta)$ and $K_{4}(\delta)$ to be the same
in the case $\sigma\geq 1$ and in the case $1/2\leq\sigma\leq 1$, it
is enough to take for each of them the maximum of the values obtained
in the two ranges.\qed \medskip

Our second a priori estimate concerns the subcritical case
$\sigma\in(0,1/2]$.  We state the result under the assumptions we need
in the sequel, even if slightly different from those of
Theorem~\ref{thm:sub-K}.

\begin{prop}[Subcritical dissipation -- A priori 
	bound]\label{prop:ap-sub}
	Let us consider equation (\ref{pbm:k-eqn}) under the following 
	assumptions:
	\begin{itemize}
		\item $A$ is a self-adjoint nonnegative operator on a
		separable Hilbert space $H$,
	
		\item $m:[0,+\infty)\to[0,+\infty)$ satisfies the strict 
		hyperbolicity assumption~(\ref{hp:k-sh}),
		
		\item $\sigma\in(0,1/2]$ and $\delta>0$ are two real 
		numbers.
	\end{itemize}
	
	Let us assume that there exists a local solution $u(t)$, defined
	in the time-interval $[0,T]$ for some $T>0$, with the regularity
	prescribed by (\ref{th:u-reg}).  Let us assume in addition that
	the continuity modulus of the function $t\to m(|A^{1/2}u(t)|^{2})$
	satisfies the assumptions of Theorem~\ref{thmbibl:sub}, and in
	particular inequality (\ref{hp:sub}) with $\Lambda_{\infty}$
	defined by (\ref{defn:Linfty}).
	
	Then it turns out that
	\begin{equation}
		|u'(t)|^{2}+|A^{1/2}u(t)|^{2}\leq
		\max\left\{1,\frac{1}{\mu_{1}}\right\}
		\left(|u_{1}|^{2}+M\left(|A^{1/2}u_{0}|^{2}\right)\right)
		\quad\quad
		\forall t\in[0,T].
		\label{th:ap-sub-k}
	\end{equation}
\end{prop}

\paragraph{\textmd{\textit{Proof}}}

Let us consider $u(t)$ as a solution to the linear equation
(\ref{pbm:lin-c(t)}) with the coefficient $c(t)$ defined as in
(\ref{defn:c(t)}).  This coefficient is bounded from below by
$\mu_{1}$ because of the strict hyperbolicity assumption, and bounded
from above by some $\mu_{2}$ because of its continuity.  Moreover, we
assumed that the continuity modulus of $c(t)$ satisfies
(\ref{hp:sub}).  As a consequence, we are in the assumptions of
Theorem~\ref{thmbibl:sub} and Theorem~\ref{thmbibl:sub-gevrey}, from
which we deduce further regularity of $u'(t)$ in $(0,T]$.

This is enough to justify rigorously the computation of the
time-derivative of $H(t)$ leading to (\ref{deriv:H'}), which therefore
holds true for every $t\in(0,T]$.  Since $H(t)$ is continuous,
equality (\ref{est:ap-easy}) holds true for every $t\in[0,T]$.

Finally, the strict hyperbolicity assumption (\ref{hp:k-sh}) implies
that $M(x)\geq\mu_{1}x$ for every $x\geq 0$.  At this point,
(\ref{th:ap-sub-k}) follows from (\ref{est:ap-easy}).\qed

\setcounter{equation}{0}
\section{Proof of main results}\label{sec:proofs-main}

We are now ready to prove Theorem~\ref{thm:sup-K} and
Theorem~\ref{thm:sub-K}.  Actually, we prove two stronger results,
which we state below.

The first one is a stronger version of Theorem~\ref{thm:sup-K} in 
interpolation spaces.

\begin{thm}[Supercritical dissipation]\label{thm:sup-K-true}
	Let us consider problem (\ref{pbm:k-eqn})--(\ref{pbm:k-data})
	under the same assumptions of Theorem~\ref{thm:sup-K}.
	
	Let us assume in addition that $(u_{0},u_{1})\in V_{\omega_{d}}$
	for some continuity modulus $\omega_{d}(x)$.

	Then the problem admits at least one global solution
	\begin{equation}
		u\in C^{0}\left([0,+\infty),D(A^{1/2})\right)\cap
		C^{1}\left([0,+\infty),H\right)
		\label{th:k-sup-true-1}
	\end{equation}
	such that
	\begin{equation}
		(u,u')\in C^{0}\left([0,+\infty),V_{\omega_{d}}\right).
		\label{th:k-sup-true-2}
	\end{equation}
\end{thm}

From Proposition~\ref{prop:data-omega}, applied with $\alpha=0$, we
know that every pair of initial conditions $(u_{0},u_{1})\in
D(A^{1/2})\times H$ lies in $V_{\omega_{d}}$ for a suitable continuity
modulus $\omega_{d}(x)$ depending on $(u_{0},u_{1})$.  As a
consequence, Theorem~\ref{thm:sup-K-true} above implies
Theorem~\ref{thm:sup-K}.
\medskip

The second result is a stronger version of Theorem~\ref{thm:sub-K}.

\begin{thm}[Subcritical dissipation]\label{thm:sub-K-true}
	Let us consider problem (\ref{pbm:k-eqn})--(\ref{pbm:k-data})
	under the same assumptions of Theorem~\ref{thm:sub-K}. 
	
	Let us assume in addition that $(u_{0},u_{1})\in V_{\omega_{d}}$
	for some continuity modulus $\omega_{d}(x)$ satisfying
	(\ref{th:od=1}) through (\ref{th:od-lim}) with the same
	$\alpha\in[0,1/4)$ which appears in (\ref{hp:limsup-alpha}).
	
	Then the following conclusions hold true.
	\begin{enumerate}
		\renewcommand{\labelenumi}{(\arabic{enumi})}
		\item  \emph{(Local existence)} There exists $T>0$ such that 
		problem (\ref{pbm:k-eqn})--(\ref{pbm:k-data}) admits at least 
		one local solution
		\begin{equation}
			u\in C^{0}\left([0,T],D(A^{1/2})\right)\cap
			C^{1}\left([0,T],H\right)
			\label{th:k-sub-true-loc1}
		\end{equation}
		such that
		\begin{equation}
			(u,u')\in C^{0}\left([0,T],V_{\omega_{d}}\right).
			\label{th:k-sub-true-loc2}
		\end{equation}
	
		\item \emph{(Alternative)} Every local solution satisfying
		(\ref{th:k-sub-true-loc1}) and (\ref{th:k-sub-true-loc2}) can
		be continued to a solution with the same regularity defined on
		a maximal time-interval $[0,T_{*})$, and either
		$T_{*}=+\infty$ or
		\begin{equation}
			\limsup_{t\to T_{*}^{-}}\|(u(t),u'(t))\|_{\omega_{d}}=+\infty.
			\label{th:limsup}
		\end{equation}
		
		\item \emph{(Global existence)} There exists $\ep_{1}>0$,
		independent of $\omega_{d}(x)$ (provided of course that
		(\ref{th:od=1}) through (\ref{th:od-lim}) are satisfied),
		with the following property.  If the initial
		conditions $(u_{0},u_{1})$ satisfy
		\begin{equation}
			|u_{1}|^{2}+|A^{1/2}u_{0}|^{2}+
			\|(u_{0},u_{1})\|_{\omega_{d}}^{2}
			\leq\ep_{1},
			\label{hp:true-data-small}
		\end{equation}
		then every local solution satisfying
		(\ref{th:k-sub-true-loc1}) and (\ref{th:k-sub-true-loc2}) can
		be continued to a global solution with the regularity
		(\ref{th:k-sup-true-1}) and (\ref{th:k-sup-true-2}).
	\end{enumerate}
\end{thm}

From Proposition~\ref{prop:data-omega}, applied with the same value of
$\alpha$ for which (\ref{hp:limsup-alpha}) holds true, we know that
every pair of initial conditions $(u_{0},u_{1})\in
D(A^{\alpha+1/2})\times D(A^{\alpha})$ lies in $V_{\omega_{d}}$ for a
suitable continuity modulus $\omega_{d}(x)$ satisfying (\ref{th:od=1})
through (\ref{th:od-data}).  As a consequence, the local existence
statement of Theorem~\ref{thm:sub-K-true} above implies the local
existence statement of Theorem~\ref{thm:sub-K}. Up to now, we did not 
exploit (\ref{th:od-data}), which only comes into play in the global 
existence statement.

If in addition the pair $(u_{0},u_{1})$ is small in the sense of
(\ref{hp:data-small}), then from (\ref{th:od-data}) it follows that
$(u_{0},u_{1})$ is small also in the sense of
(\ref{hp:true-data-small}) with $\ep_{1}=2\ep_{0}$.
Since $\ep_{1}$ does not depend on $\omega_{d}(x)$, the global
existence statement of Theorem~\ref{thm:sub-K-true} above implies the
global existence statement of Theorem~\ref{thm:sub-K}.

In conclusion, our main results are true if we prove 
Theorem~\ref{thm:sup-K-true} and Theorem~\ref{thm:sub-K-true}. 

\subsection*{Proof of Theorem~\ref{thm:sup-K-true}}

\paragraph{\textmd{\textit{Local existence}}}

Let us set
\begin{equation}
	M_{0}:=|A^{1/2}u_{0}|^{2}+1,
	\quad\quad
	\mu_{2}:=\max_{x\in[0,M_{0}]}m(x),
	\quad\quad
	m_{*}(x):=m\left(\min\{x,M_{0}\}\strut\right).
	\label{defn:true-const}
\end{equation}

In such a way, $m_{*}(x)$ is a nonnegative continuous function, which
is bounded from above by $\mu_{2}$, and coincides with $m(x)$ for
every $x\in[0,M_{0}]$.  Moreover, the continuity modulus of $m_{*}(x)$
in $[0,+\infty)$ coincides with the continuity modulus of $m(x)$ in
$[0,M_{0}]$, which we denote by $\omega_{m}(x)$.

We claim that problem (\ref{pbm:k-eqn})--(\ref{pbm:k-data}), with
$m_{*}(x)$ instead of $m(x)$, admits at least one local solution
$u_{*}(t)$, defined in some time-interval $[0,T]$, with regularity
\begin{equation}
	u_{*}\in C^{0}\left([0,T],D(A^{1/2})\right)\cap
	C^{1}\left([0,T],H\right),
	\label{th:k-sup-true-loc1}
\end{equation}
and
\begin{equation}
	(u_{*},u_{*}')\in C^{0}\left([0,T],V_{\omega_{d}}\right).
	\label{th:k-sup-true-loc2}
\end{equation}

If we prove this claim, then $u_{*}(t)$ turns out to be also a
solution to the original problem in a possibly smaller time-interval.
Indeed, since $|A^{1/2}u_{*}(t)|^{2}$ is a continuous function and
$|A^{1/2}u_{*}(0)|^{2}=|A^{1/2}u_{0}|^{2}<M_{0}$, there exists
$T_{1}\in(0,T]$ such that $|A^{1/2}u_{*}(t)|^{2}\leq M_{0}$ for every
$t\in[0,T_{1}]$.  As a consequence 
$$m_{*}(|A^{1/2}u_{*}(t)|^{2})=
m(|A^{1/2}u_{*}(t)|^{2}) 
\quad\quad
\forall t\in[0,T_{1}],$$
and hence $u_{*}(t)$ is also a solution to the original problem in
$[0,T_{1}]$.

In order to prove the existence of $u_{*}(t)$, we exploit a fixed
point argument.  Let us choose a positive number $T$, and a real
number $\nu\geq 1$ satisfying (\ref{hp:sup-nu}).   Let us set
$$\omega_{c}(x):=\omega_{m}\left(3K_{1}(\delta,\mu_{2})
e^{(1+\mu_{2})\nu T}\cdot
\|(u_{0},u_{1})\|_{\omega_{d}}^{2}\cdot
\omega_{d}(x)\strut\right),$$
where $K_{1}(\delta,\mu_{2})$ is the constant of
Theorem~\ref{thmbibl:sup}.  Let $\mathbb{X}$ be the space of all
continuous functions $c:[0,T]\to[0,\mu_{2}]$ with continuity modulus
equal to $\omega_{c}(x)$.  Due to Ascoli's Theorem, the space
$\mathbb{X}$ turns out to be a compact and convex subset of the Banach
space of all continuous functions $c:[0,T]\to\re$, endowed with the
sup norm.

For every $c\in \mathbb{X}$, we consider the solution $u(t)$ to the
linear problem (\ref{pbm:lin-c(t)})--(\ref{pbm:k-data}), and then we
define a function $\Phi(c):[0,T]\to\re$ as
\begin{equation}
	\left[\Phi(c)\right](t):=m_{*}(|A^{1/2}u(t)|^{2})
	\quad\quad
	\forall t\in[0,T].
	\label{defn:Phi(c)}
\end{equation}

We claim that $\Phi$ is a continuous map from $\mathbb{X}$ to
$\mathbb{X}$.  If we prove this claim, then Schauder's fixed point
theorem implies that $\Phi$ admits at least one fixed point $c(t)$.
The corresponding solution to the linear problem is the solution
$u_{*}(t)$ we were looking for, and its regularity
(\ref{th:k-sup-true-1}) follows from Theorem~\ref{thmbibl:sup}, while
its regularity (\ref{th:k-sup-true-2}) follows from
Proposition~\ref{prop:sup-sol-omega}.

So let us verify the properties of $\Phi$.  First of all, we need to
show that $\Phi(c)\in\mathbb{X}$ for every $c\in\mathbb{X}$.  The fact
that $[\Phi(c)](t)\in[0,\mu_{2}]$ for every $t\in[0,T]$ follows from
the bounds on $m_{*}(x)$.  The main point is proving that the map
$t\to [\Phi(c)](t)$ has continuity modulus equal to $\omega_{c}(x)$.
To this end, we begin by applying
Proposition~\ref{prop:sup-sol-omega}, from which we deduce that
$(u(t),u'(t))\in V_{\omega_{d}}$ for every $t\in[0,T]$ and
(\ref{th:od-cont}) holds true.  This allows us to apply
Proposition~\ref{prop:o-cont-auq}, from which we deduce that
$$\left||A^{1/2}u(a)|^{2}-|A^{1/2}u(b)|^{2}\right|\leq
3K_{1}(\delta,\mu_{2})e^{(1+\mu_{2})\nu T}\cdot
\|(u_{0},u_{1})\|_{\omega_{d}}^{2}\cdot
\omega_{d}(|a-b|)$$
for every $a$ and $b$ in $[0,T]$. Since  $m_{*}(x)$ has 
continuity modulus $\omega_{m}(x)$, we conclude that
\begin{eqnarray*}
	\left|\strut\left[\Phi(c)\right](a)-\left[\Phi(c)\right](b)\right| & 
	= & \left|m_{*}\left(|A^{1/2}u(a)|^{2}\right)-
	m_{*}\left(|A^{1/2}u(b)|^{2}\right)\right|  \\
	 & \leq & \omega_{m}\left(
	 \left||A^{1/2}u(a)|^{2}-|A^{1/2}u(b)|^{2}\right|\right)  \\
	 & \leq & \omega_{m}\left(\strut
	 3K_{1}(\delta,\mu_{2})e^{(1+\mu_{2})\nu T}
	 \cdot\|(u_{0},u_{1})\|_{\omega_{d}}^{2}\cdot
	 \omega_{d}(|a-b|)\right)  \\
	 & = & \omega_{c}(|a-b|)
\end{eqnarray*}
for every $a$ and $b$ in $[0,T]$, which completes the proof that
$\Phi(c)\in\mathbb{X}$.

Now that we know that $\Phi$ maps $\mathbb{X}$ to $\mathbb{X}$, its 
continuity with respect to the sup norm is a simple application of 
Theorem~\ref{thm:sup-cn2cinfty}. This completes the proof of local 
existence.

\paragraph{\textmd{\textit{Alternative}}}

Every local solution to (\ref{pbm:k-eqn})--(\ref{pbm:k-data}) 
satisfying (\ref{th:k-sup-true-loc1}) and (\ref{th:k-sup-true-loc2}) can be 
continued to a maximal solution $u(t)$ defined in some time-interval 
$[0,T_{*})$, and such that 
\begin{equation}
	u\in C^{0}\left([0,T_{*}),D(A^{1/2})\right)\cap
	C^{1}\left([0,T_{*}),H\right),
	\label{th:k-sub-true-*1}
\end{equation}
\begin{equation}
	(u,u')\in C^{0}\left([0,T_{*}),V_{\omega_{d}}\right).
	\label{th:k-sub-true-*2}
\end{equation}

We prove that either $T_{*}=+\infty$ or
\begin{equation}
	\limsup_{t\to T_{*}^{-}}\|(u(t),u'(t))\|_{\omega_{d}}^{2}=+\infty.
	\label{th:sup-limsup}
\end{equation}

To this end, let us assume by contradiction that 
$T_{*}<+\infty$ and (\ref{th:sup-limsup}) is false, so that
\begin{equation}
	L_{0}:=\sup_{t\in[0,T_{*})}\|(u(t),u'(t))\|_{\omega_{d}}^{2}<+\infty.
	\label{limsup-false}
\end{equation}

We claim that $u(t)$ can be extended to the closed interval 
$[0,T_{*}]$ in such a way that
\begin{equation}
	u\in C^{0}\left([0,T_{*}],D(A^{1/2})\right)\cap
	C^{1}\left([0,T_{*}],H\right),
	\label{th:k-sub-true-*1]}
\end{equation}
with in addition
\begin{equation}
	(u,u')\in C^{0}\left([0,T_{*}],V_{\omega_{d}}\right).
	\label{th:k-sub-true-*2]}
\end{equation}

If we prove this claim, then we can use $(u(T_{*}),u'(T_{*}))$ as a
new initial condition and extend the solution beyond $T_{*}$, thus
contradicting its maximality.

In order to prove the claim, we apply
Proposition~\ref{prop:o-cont-auq} in every closed interval
$[0,T]\subseteq[0,T_{*})$.  We deduce that the function
$t\to|A^{1/2}u(t)|^{2}$ is continuous in $[0,T]$ with continuity
modulus equal to $3L_{0}\omega_{d}(x)$.  Since the continuity modulus
does not depend on $T$, the function is continuous with the same
continuity modulus in $[0,T_{*})$, and in particular it is bounded in
the same interval.  As a consequence, the function $t\to
c(t):=m(|A^{1/2}u(t)|^{2})$ is uniformly continuous in $[0,T_{*})$,
and hence it can be extended to a continuous function defined in the
closed interval $[0,T_{*}]$.

The function $c(t)$, in the closed interval $[0,T_{*}]$, is
nonnegative and bounded from above.  Therefore, from
Theorem~\ref{thmbibl:sup} and Proposition~\ref{prop:sup-sol-omega} it
follows that the linear problem
(\ref{pbm:lin-c(t)})--(\ref{pbm:k-data}) has a unique solution in the
closed interval $[0,T_{*}]$ with the regularity
(\ref{th:k-sub-true-*1]}) and (\ref{th:k-sub-true-*2]}), which is the
required extension of $u$.

\paragraph{\textmd{\textit{Global existence}}}

Let $u(t)$ be a local solution defined in a maximal interval 
$[0,T_{*})$. Let us assume by contradiction that $T_{*}<+\infty$. Due 
to the a priori estimate of Proposition~\ref{prop:ap-sup}, we know that 
$|A^{1/2}u(t)|$ is uniformly bounded in $[0,T_{*})$, hence also 
$c(t)$ is uniformly bounded from above in $[0,T_{*})$ by some 
constant $\widehat{\mu}_{2}$. Therefore, from 
Proposition~\ref{prop:sup-sol-omega} we deduce that
$$\|(u(t),u'(t))\|_{\omega_{d}}^{2}\leq
K_{1}(\delta,\widehat{\mu}_{2})e^{(1+\widehat{\mu}_{2})\nu t}
\cdot\|(u_{0},u_{1})\|_{\omega_{d}}^{2}
\quad\quad
\forall t\in[0,T_{*}).$$

Passing to the limit as $t\to T_{*}^{-}$ we contradict 
(\ref{th:sup-limsup}).\qed

\subsection*{Proof of Theorem~\ref{thm:sub-K-true}}

\paragraph{\textmd{\textit{Local existence}}}

The argument is similar to the proof of Theorem~\ref{thm:sup-K-true}. 
Once again, we define $M_{0}$, $\mu_{2}$, $m_{*}(x)$ as in 
(\ref{defn:true-const}), and we look for a solution $u_{*}(t)$ to the 
modified problem, which once again turns out to be a solution to the 
original problem in a possibly smaller time-interval. We remark that 
in this case $m_{*}(x)$ is bounded from below by $\mu_{1}$ due to 
the strict hyperbolicity assumption (\ref{hp:k-sh}).

The existence of $u_{*}(t)$ follows again from Schauder's fixed point 
theorem. This time we set
$$\omega_{c}(x):=\omega_{m}\left(6K_{2}(\delta,\mu_{1},\mu_{2})
\cdot\|(u_{0},u_{1})\|_{\omega_{d}}^{2}
\cdot\omega_{d}(x)\strut\right),$$
where $K_{2}(\delta,\mu_{1},\mu_{2})$ is the constant of 
Theorem~\ref{thmbibl:sub}. First of all, we prove that
\begin{equation}
	\lim_{x\to 0^{+}}\frac{\omega_{m}(\omega_{d}(x))}{x^{1-2\sigma}}=0.
	\label{th:omega-md}
\end{equation}

This is trivial if $\sigma=1/2$.  If $\sigma>1/2$, the value $\alpha$
for which (\ref{hp:limsup-alpha}) holds true is positive.  Therefore,
we can write
$$\frac{\omega_{m}(\omega_{d}(x))}{x^{1-2\sigma}}=
\left(\frac{[\omega_{m}(\omega_{d}(x))]^{4\alpha}}
{[\omega_{d}(x)]^{1-2\sigma}}\right)^{1/(4\alpha)}
\left(\frac{\omega_{d}(x)}{x^{4\alpha}}\right)^{(1-2\sigma)/(4\alpha)},$$
and observe that the first term is bounded because of
(\ref{hp:limsup-alpha}), while the second term tends to 0 as $x\to
0^{+}$ because of (\ref{th:od-lim}).

From (\ref{th:omega-md}) and (\ref{hp:omega-L}) it follows that
\begin{equation}
	\lim_{x\to 0^{+}}\frac{\omega_{c}(x)}{x^{1-2\sigma}}=0,
	\label{omega-c-lim}
\end{equation}
and hence there exists $\nu>0$ such that
\begin{equation}
	4\delta^{2}\mu_{1}\geq
	\left[\lambda^{1-2\sigma}\omega_{c}\left(\frac{1}{\lambda}\right)\right]^{2}+
	2\delta\left[\lambda^{1-2\sigma}\omega_{c}\left(\frac{1}{\lambda}\right)\right]
	\quad\quad
	\forall\lambda\geq\nu.
	\label{hp:sub-nuc}
\end{equation}

Given this value of $\nu$, we finally set
\begin{equation}
	T:=\frac{\log 2}{(1+\mu_{2})\nu}.
	\label{defn:T}
\end{equation}

Now we proceed again in analogy with the proof of
Theorem~\ref{thm:sup-K-true}.  We consider the space $\mathbb{X}$ of
all functions $c:[0,T]\to[\mu_{1},\mu_{2}]$ with continuity modulus
$\omega_{c}(x)$, and for every $c\in\mathbb{X}$ we define
$\Phi(c):[0,T]\to\re$ as in (\ref{defn:Phi(c)}).  The solution $u(t)$
to the linear problem (\ref{pbm:lin-c(t)})--(\ref{pbm:k-data}) needed
in the definition is now provided by Theorem~\ref{thmbibl:sub}, whose
assumptions are satisfied because $\omega_{c}(x)$ satisfies
(\ref{omega-c-lim}).

Once again, the main point is proving that the continuity modulus of
the map $t\to [\Phi(c)](t)$ is $\omega_{c}(x)$.  To this end, we begin
by applying Proposition~\ref{prop:sub-sol-omega}, from which we deduce
that $(u,u')\in C^{0}([0,T],V_{\omega_{d}})$ and
(\ref{th:od-sub-cont}) holds true.  Therefore, from (\ref{defn:T}) it
follows that
$$\|(u(t),u'(t))\|_{\omega_{d}}^{2}\leq
2K_{2}(\delta,\mu_{1},\mu_{2})\|(u_{0},u_{1})\|_{\omega_{d}}^{2}
\quad\quad
\forall t\in[0,T].$$

This allows us to apply Proposition~\ref{prop:o-cont-auq}, from which 
we deduce that
$$\left||A^{1/2}u(a)|^{2}-|A^{1/2}u(b)|^{2}\right|\leq
6K_{2}(\delta,\mu_{1},\mu_{2})\|(u_{0},u_{1})\|_{\omega_{d}}^{2}\cdot
\omega_{d}(|a-b|)$$
for every $a$ and $b$ in $[0,T]$.  At this point we conclude 
as in the proof of Theorem~\ref{thm:sup-K-true}.

\paragraph{\textmd{\textit{Alternative}}}

We argue as in the corresponding paragraph of the proof of
Theorem~\ref{thm:sup-K-true}.  Assuming by contradiction that
$T_{*}<+\infty$ and (\ref{limsup-false}) holds true, the key point is
extending $c(t)$ to the closed interval $[0,T_{*}]$ in such a way that
the assumptions of Theorem~\ref{thmbibl:sub} and
Proposition~\ref{prop:sub-sol-omega} are satisfied, so that also the
solution $u(t)$ can be continued to $[0,T_{*}]$ with the required
regularity.

In analogy with the proof of Theorem~\ref{thm:sup-K-true}, from
Proposition~\ref{prop:o-cont-auq} we deduce that the function
$t\to|A^{1/2}u(t)|^{2}$ is uniformly continuous in $[0,T_{*})$, with
continuity modulus equal to $3L_{0}\omega_{d}(x)$.  As a consequence,
the function $t\to c(t)$ is uniformly continuous in $[0,T_{*})$, with
continuity modulus $\omega_{m}(3L_{0}\omega_{d}(x))$, and hence it can
be extended to the closed interval $[0,T_{*}]$ with the same
continuity modulus.  On the other hand, this continuity modulus
satisfies (\ref{hp:sub}) with $\Lambda_{\infty}=0$ (same proof as in
the case of $\omega_{c}(x)$, just with a different constant).
Therefore, the assumptions of Theorem~\ref{thmbibl:sub} are satisfied.

\paragraph{\textmd{\textit{Global existence}}}

Let us consider a new continuity modulus defined by
$$\widehat{\omega}_{c}(x):=\omega_{m}(6L_{1}x^{4\alpha}),$$
where $L_{1}>0$ is a parameter to be chosen. We claim that
$\widehat{\omega}_{c}(x)$ satisfies (\ref{hp:sub-nuc}) with $\nu=1$, 
namely
\begin{equation}
	4\delta^{2}\mu_{1}\geq
	\left[\lambda^{1-2\sigma}\widehat{\omega}_{c}
	\left(\frac{1}{\lambda}\right)\right]^{2}+
	2\delta\left[\lambda^{1-2\sigma}\widehat{\omega}_{c}
	\left(\frac{1}{\lambda}\right)\right]
	\quad\quad
	\forall\lambda\geq 1,
	\label{hp:sub-nuc-hat}
\end{equation}
provided that $L_{1}$ is small enough.  If $\sigma=1/2$, this is
almost obvious.  If $\sigma<1/2$, the value $\alpha$ for which
(\ref{hp:limsup-alpha}) holds true is positive, and hence
$$R:=\sup_{y\in(0,1]}\frac{\omega_{m}(y^{4\alpha})}{y^{1-2\sigma}}=
\sup_{z\in(0,1]}
\left(\frac{[\omega_{m}(z)]^{4\alpha}}{z^{1-2\sigma}}\right)^{1/4\alpha}
<+\infty.$$

Therefore, when $6L_{1}\leq 1$ it follows that
$$\sup_{x\in(0,1]}\frac{\widehat{\omega}_{c}(x)}{x^{1-2\sigma}}=
\sup_{x\in(0,1]}
\frac{\omega_{m}\left([(6L_{1})^{1/(4\alpha)}x]^{4\alpha}\right)}
{[(6L_{1})^{1/(4\alpha)}x]^{1-2\sigma}}\cdot
(6L_{1})^{(1-2\sigma)/(4\alpha)}\leq
R\cdot(6L_{1})^{(1-2\sigma)/(4\alpha)},$$
which implies (\ref{hp:sub-nuc-hat}) when $L_{1}$ is small enough.

Now let us choose $L_{1}>0$ so that (\ref{hp:sub-nuc-hat}) holds true,
let us set
$$\widehat{\mu}_{2}:=\max\{m(x):x\in[0,L_{1}]\},$$
and let us assume that $\ep_{1}$ in (\ref{hp:true-data-small}) is 
small enough so that
\begin{equation}
	K_{2}(\delta,\mu_{1},\widehat{\mu}_{2})
	\|(u_{0},u_{1})\|_{\omega_{d}}^{2}+
	\max\left\{1,\frac{1}{\mu_{1}}\right\}\left(
	|u_{1}|^{2}+M\left(|A^{1/2}u_{0}|^{2}\right)
	\right)\leq L_{1},
	\label{hp:small-L1}
\end{equation}
where $K_{2}(\delta,\mu_{1},\widehat{\mu}_{2})$ is once again the
constant of Theorem~\ref{thmbibl:sub}. We point out that this 
smallness condition on $\ep_{1}$ is independent of $\omega_{d}(x)$.

For any such initial condition, let us consider a local solution $u(t)$ to 
(\ref{pbm:k-eqn})--(\ref{pbm:k-data}), defined on a maximal interval 
$[0,T_{*})$, with the regularity prescribed by (\ref{th:k-sub-true-*1}) and 
(\ref{th:k-sub-true-*2}). 

We claim that
\begin{equation}
	\|(u(t),u'(t))\|_{\omega_{d}}^{2}\leq 2L_{1}
	\quad\quad
	\forall t\in[0,T_{*}).
	\label{est:2L}
\end{equation}

Since this contradicts (\ref{th:limsup}), the only possibility left by
the alternative is that $T_{*}=+\infty$, which means that $u(t)$ is a 
global solution.  In order to prove (\ref{est:2L}), we introduce
$$S:=\sup\left\{t\in[0,T_{*}): \|(u(\tau),u'(\tau))\|_{\omega_{d}}\leq
2L_{1}\quad\forall\tau\in[0,t]\strut\right\},$$
so that (\ref{est:2L}) is now equivalent to proving that $S=T_{*}$. 

From (\ref{hp:small-L1}) it follows that
$\|(u(0),u'(0))\|_{\omega_{d}}\leq L_{1}<2L_{1}$ (because 
$K_{2}(\delta,\mu_{1},\widehat{\mu}_{2})\geq 1$, as already observed after 
Theorem~\ref{thmbibl:sub}).  Therefore, $S$ is
the supremum of a nonempty set, and $S>0$ because the function
$t\to\|(u(t),u'(t))\|_{\omega_{d}}$ is continuous.  Moreover, the
maximality of $S$ implies that either $S=T_{*}$, or
\begin{equation}
	\|(u(S),u'(S))\|_{\omega_{d}}^{2}= 2L_{1}.
	\label{est:=2L}
\end{equation}

If we prove that (\ref{est:=2L}) does not hold, then automatically
$S=T_{*}$.  So let us assume by contradiction that $S<T_{*}$.  From
our definition of $S$ it follows that
$$\|(u(t),u'(t))\|_{\omega_{d}}^{2}\leq 2L_{1}
\quad\quad
\forall t\in[0,S].$$

Therefore, from Proposition~\ref{prop:o-cont-auq} we obtain that the
function $t\to|A^{1/2}u(t)|^{2}$ has continuity modulus
$6L_{1}\omega_{d}(x)$ in $[0,S]$.  As a consequence, the function
\begin{equation}
	c(t):=m\left(|A^{1/2}u(t)|^{2}\right)
	\label{defn:c(t)-bis}
\end{equation}
has continuity modulus $\omega_{m}(6L_{1}\omega_{d}(x))$ in $[0,S]$,
and this continuity modulus is less than or equal to
$\widehat{\omega}_{c}(x)$ because of (\ref{th:od<=}).  Thus the
continuity modulus of $c(t)$ satisfies (\ref{hp:sub}), and hence we
can apply Proposition~\ref{prop:ap-sub} and deduce that
(\ref{th:ap-sub-k}) holds true for every $t\in[0,S]$.  In particular,
$|A^{1/2}u(t)|^{2}\leq L_{1}$ and hence $c(t)\leq\widehat{\mu}_{2}$
for every $t\in[0,S]$.

Now we estimate the (semi)norm of $(u(t),u'(t))$ in $V_{\omega_{d}}$.
To this end, we interpret $u(t)$ as a solution to the linear equation
(\ref{pbm:lin-c(t)}) with coefficient $c(t)$ given by
(\ref{defn:c(t)-bis}).  Since $c(t)$ has continuity modulus
(less than or equal to) 
$\widehat{\omega}_{c}(x)$, the assumptions of
Proposition~\ref{prop:sub-sol-omega} are satisfied with $\nu:=1$.  As
for low-frequency components, from (\ref{th:od-sub-cont-}) with
$\nu=1$ it follows that
\begin{eqnarray*}
	\|(u_{\nu=1,-}(t),u_{\nu=1,-}'(t))\|_{\omega_{d}}^{2} & \leq & 
	\frac{1}{\omega_{d}(1)}\left(
	|u_{\nu=1,-}'(t)|^{2}+|A^{1/2}u_{\nu=1,-}(t)|^{2} \right) \\
	 & \leq & \frac{1}{\omega_{d}(1)}\left(
	|u'(t)|^{2}+|A^{1/2}u(t)|^{2} \right).
\end{eqnarray*}

Keeping into account that $\omega_{d}(1)=1$, from the a priori
estimate of Proposition~\ref{prop:ap-sub} we deduce that
\begin{equation}
	\|(u_{\nu=1,-}(t),u_{\nu=1,-}'(t))\|_{\omega_{d}}^{2}\leq
	\max\left\{1,\frac{1}{\mu_{1}}\right\}\left(
	|u_{1}|^{2}+M\left(|A^{1/2}u_{0}|^{2}\right)\right).
	\label{est:ge-u-}
\end{equation}

As for high-frequency components, from (\ref{th:od-sub-cont+}) it
follows that
\begin{equation}
	\|(u_{\nu=1,+}(t),u_{\nu=1,+}'(t))\|_{\omega_{d}}^{2} \leq
	K_{2}(\delta,\mu_{1},\widehat{\mu}_{2})
	\|(u_{0},u_{1})\|_{\omega_{d}}^{2}.
	\label{est:ge-u+}
\end{equation}

Summing (\ref{est:ge-u-}) and (\ref{est:ge-u+}), and keeping
(\ref{hp:small-L1}) into account, we finally obtain that
$\|(u(t),u'(t))\|_{\omega_{d}}^{2}\leq L_{1}$ for every $t\in[0,S]$.
This contradicts (\ref{est:=2L}).\qed

\label{NumeroPagine}


{\small 

\begin{thebibliography}{99}

	\bibitem{ap}{\sc A.\ Arosio, S.\ Panizzi}; On the well-posedness of
	the Kirchhoff string.  {\em Trans.\ Amer.\ Math.\ Soc.}\ \textbf{348}
	(1996), no.~1, 305--330.

	\bibitem{as}{\sc A.\ Arosio, S.\ Spagnolo}; Global solutions to
	the Cauchy problem for a nonlinear hyperbolic equation.
	\emph{Nonlinear partial differential equations and their
	applications.  Coll\`{e}ge de France seminar, Vol.\ VI (Paris,
	1982/1983)}, 1--26, Res.\ Notes in Math., 109, Pitman, Boston, MA,
	1984.

	\bibitem{bernstein}{\sc S.~Bernstein}; Sur une classe
	d'\'equations fonctionnelles aux d\'eriv\'ees partielles,
	(Russian, French summary) {\em Bull.\ Acad.\ Sci.\ URSS. S\'{e}r.\
	Math.\ [Izvestia Akad.\ Nauk SSSR]} \textbf{4} (1940),  17--26.

	\bibitem{CR}\textsc{G.~Chen, D.~L.~Russell}; A mathematical model
	for linear elastic systems with structural damping.  \emph{Quart.\
	Appl.\ Math.}\ \textbf{39} (1981/82), no.~4, 433--454.

	\bibitem{CT1}\textsc{S.~P.~Chen, R.~Triggiani}; Proof of
	extensions of two conjectures on structural damping for elastic
	systems.  \emph{Pacific J.\ Math.}\ \textbf{136} (1989), no.~1,
	15--55.

	\bibitem{CT2}\textsc{S.~P.~Chen, R.~Triggiani}; Characterization
	of domains of fractional powers of certain operators arising in
	elastic systems, and applications.  \emph{J.\ Differential
	Equations} \textbf{88} (1990), no.~2, 279--293.
	
	\bibitem{CT3}\textsc{S.~P.~Chen, R.~Triggiani}; Gevrey class
	semigroups arising from elastic systems with gentle dissipation:
	the case $0<\alpha<1/2$.  \emph{Proc.\ Amer.\ Math.\ Soc.}\
	\textbf{110} (1990), no.~2, 401--415.
	
	\bibitem{dgcs}{\sc F.\ Colombini, E.\ De Giorgi, S.\ Spagnolo}; Sur le
	\'{e}quations hyperboliques avec des coefficients qui ne d\'{e}pendent
	que du temp.  (French) {\em Ann.\ Scuola Norm.\ Sup.\ Pisa Cl.\ Sci.\
	(4)} \textbf{6} (1979), no.~3, 511--559.

	\bibitem{ds:analytic-1}{\sc P.\ D'Ancona, S.\ Spagnolo}; Global
	solvability for the degenerate Kirchhoff equation with real
	analytic data.  {\em Invent.\ Math.}\ \textbf{108} (1992), no.~2,
	247--262.

	\bibitem{ds:analytic-2}{\sc P.\ D'Ancona, S.\ Spagnolo}; On an abstract
	weakly hyperbolic equation modelling the nonlinear vibrating
	string.  {\em Developments in partial differential equations and
	applications to mathematical physics (Ferrara, 1991)}, 27--32,
	Plenum, New York, 1992.

	\bibitem{ds:dispersive}{\sc P.\ D'Ancona, S.\ Spagnolo}; A class of
	nonlinear hyperbolic problems with global solutions.  {\em Arch.\
	Rational Mech.\ Anal.}\ \textbf{124} (1993), no.~3, 201--219.

	\bibitem{debrito}{\sc E.\ H.\ de Brito}; The damped elastic stretched
	string equation generalized: existence, uniqueness, regularity and
	stability.  \emph{Applicable Anal.}\ \textbf{13} (1982), no.~3,
	219--233.

	\bibitem{gg:diss}{\sc M.\ Ghisi, M.\ Gobbino}; Global
	existence and asymptotic behavior for a mildly degenerate dissipative
	hyperbolic equation of Kirchhoff type.  \emph{Asymptot.\  Anal.}\ 
	\textbf{40} (2004), no.~1, 25--36.

	\bibitem{gg:der-loss}{\sc M.\ Ghisi, M.\ Gobbino}; Derivative loss
	for Kirchhoff equations with non-Lipschitz nonlinear term.
	\emph{Ann.\ Scuola Norm.\ Sup.\ Pisa Cl.\ Sci.\ (5)} \textbf{8}
	(2009), no.~4, 613--646.

	\bibitem{gg:global}{\sc M.\ Ghisi, M.\ Gobbino}; Spectral gap
	global solutions for degenerate Kirchhoff equations.
	\emph{Nonlinear Anal.}\ \textbf{71} (2009), no.~9, 4115--4124.

	\bibitem{gg:survey-diss}{\sc M.\ Ghisi, M.\ Gobbino};
	Hyperbolic-parabolic singular perturbation for Kirchhoff equations
	with weak dissipation.  \emph{Rend.\ Ist.\ Mat.\ Univ.\ Trieste}
	\textbf{42} Suppl.\ (2010), 67--88.
	
	\bibitem{gg:survey}{\sc M.\ Ghisi, M.\ Gobbino}; Kirchhoff
	equations in generalized Gevrey spaces: local existence, global
	existence, uniqueness. {\em Rend.\ Istit.\ Mat.\ Univ.\ Trieste}.
	\textbf{42} Suppl.\ (2010), 89--110.

	\bibitem{gg:qa}{\sc M.\ Ghisi, M.\ Gobbino}; Kirchhoff equations
	from quasi-analytic to spectral-gap data.  \emph{Bull.\ Lond.\
	Math.\ Soc.}\ \textbf{43} (2011), no.~2, 374--385.

	\bibitem{gg:dgcs-strong}{\sc M.\ Ghisi, M.\ Gobbino};
	Linear hyperbolic equations with time-dependent propagation
	speed and strong damping. Preprint 	\texttt{arXiv:1408.3499 [math.AP]}.

	\bibitem{ggh:sd}{\sc M.\ Ghisi, M.\ Gobbino, H.\ Haraux}; Local
	and global smoothing effects for some linear hyperbolic equations
	with a strong dissipation.  \emph{Trans.\ Amer.\ Math.\ Soc.}\ To
	appear.  Preprint \texttt{arXiv:1402.6595 [math.AP]}.

	\bibitem{gh}{\sc J.\ M.\ Greenberg, S.\ C.\ Hu}; The initial value
	problem for a stretched string.  {\em Quart.\ Appl.\ Math.}\
	\textbf{38} (1980/81), no.~3, 289--311.

	\bibitem{HO}\textsc{A.~Haraux, M.~\^{O}tani}; Analyticity and
	regularity for a class of second order evolution equation.
	\emph{Evol.\ Equat.\ Contr.\ Theor.}\ \textbf{2} (2013),
	no.~1, 101--117.

	\bibitem{hirosawa} {\sc F.\ Hirosawa}; Global solvability for
	Kirchhoff equation in special classes of non-analytic functions.
	\emph{J.\ Differential Equations} \textbf{230} (2006), no.~1,
	49--70.

	\bibitem{I2}\textsc{R.~Ikehata, M.~Natsume}; Energy decay estimates
	for wave equations with a fractional damping.  \emph{ Differential
	Integral Equations} \textbf{25} (2012), no.~9-10, 939--956.

	\bibitem{I3}\textsc{R.~Ikehata, G.~Todorova, B.~Yordanov}; Wave
	equations with strong damping in Hilbert spaces.  \emph{J.\
	Differential Equations} \textbf{254} (2013), no.~8, 3352--3368.

	\bibitem{kirchhoff} {\sc G.\ Kirchhoff}; \emph{Vorlesungen ober
	mathematische Physik: Mechanik} (section 29.7), Teubner, Leipzig,
	1876.
	
	\bibitem{manfrin}{\sc R.\ Manfrin}; On the global solvability of
	Kirchhoff equation for non-analytic initial data.  {\em J.\
	Differential Equations} \textbf{211} (2005), no.~1, 38--60.

	\bibitem{matos-pereira}\textsc{M.~P.~Matos, D.~C.~Pereira}; On a
	hyperbolic equation with strong damping.  \emph{Funkcial.\
	Ekvac.}\ \textbf{34} (1991), no.~2, 303--311.
	
	\bibitem{matsu-ruzha}\textsc{T.~Matsuyama, M.~Ruzhansky}; Global
	well-posedness of Kirchhoff systems.  \emph{J.\ Math.\ Pures
	Appl.\ (9)} \textbf{100} (2013), no.~2, 220--240.
	
	\bibitem{nishihara}{\sc K. Nishihara}; On a global solution of
	some quasilinear hyperbolic equation.  \emph{Tokyo J.\ Math.}\
	\textbf{7} (1984), no.~2, 437--459.

	\bibitem{nishihara:k-strong}\textsc{K.~Nishihara}; Degenerate
	quasilinear hyperbolic equation with strong damping.
	\emph{Funkcial.\ Ekvac.}\ \textbf{27} (1984), no.~1, 125--145.

	\bibitem{nishihara:decay}\textsc{K.~Nishihara}; Decay properties of
	solutions of some quasilinear hyperbolic equations with strong
	damping.  \emph{Nonlinear Anal.}\ \textbf{21} (1993), no.~1,
	17--21.

	\bibitem{ny}{\sc K.\ Nishihara, Y.\ Yamada}; On global solutions of
	some degenerate quasilinear hyperbolic equations with dissipative
	terms.  \emph{Funkcial.\  Ekvac.}\ \textbf{33} (1990), no.~1, 151--159.

	\bibitem{ono}\textsc{K.~Ono}; On decay properties of solutions for
	degenerate strongly damped wave equations of Kirchhoff type.
	\emph{J.\ Math.\ Anal.\ Appl.}\ \textbf{381} (2011), no.~1,
	229--239.
	
	\bibitem{poho-m}{\sc S.\ I.\ Pohozaev}; The Kirchhoff quasilinear
	hyperbolic equation, \emph{Differentsial'nye Uravneniya}
	\textbf{21} (1985),  101--108 (English transl.:
	\emph{Differential Equations} \textbf{21} (1985), 82--87).

	\bibitem{reed}{\sc M.\ Reed, B.\ Simon}; \emph{Methods of Modern
	Mathematical Physics, I: Functional Analysis.  Second edition}.
	Academic Press, New York, 1980.
	
	\bibitem{shibata}\textsc{Y.~Shibata}; On the rate of decay of
	solutions to linear viscoelastic equation.  \emph{Math.\ Methods
	Appl.\ Sci.}\ \textbf{23} (2000), no.~3, 203--226.

	\bibitem{yamada}{\sc Y.\ Yamada}; On some quasilinear wave equations with
	dissipative terms.  \emph{Nagoya Math.\ J.}\ \textbf{87} (1982),
	17--39.

	\bibitem{yamazaki:rdg}{\sc T.\ Yamazaki}; On local solutions of some
	quasilinear degenerate hyperbolic equations. \emph{Funkcial.\
	Ekvac.}\ \textbf{31} (1988), no.~3, 439--457.

	\bibitem{yamazaki:external}{\sc T.\ Yamazaki}; Global solvability for the
	Kirchhoff equations in exterior domains of dimension three.
	\emph{J.\ Differential Equations} \textbf{210} (2005), no.~2,
	290--316.

\end{thebibliography}
}
\end{document}